\documentclass[a4paper,french,12pt]{amsart}
\usepackage{bourbaki}
\usepackage[all]{xy}
\usepackage[frenchb]{babel}
\usepackage[applemac]{inputenc}
\usepackage{eucal}
\usepackage{graphicx}
\usepackage{epsfig}
\usepackage{mathrsfs}
\usepackage{amssymb}
\usepackage{amsxtra}
\usepackage{enumerate}
\newcommand{\deactivateaddcontentsline}{%
  \let\addcontentslineoriginal\addcontentsline
  \renewcommand\addcontentsline[3]{}%
}
\newcommand{\reactivateaddcontentsline}{%
  \let\addcontentsline\addcontentslineoriginal
}

\makeatletter
\let\@@seccntformat\@seccntformat
\renewcommand*{\@seccntformat}[1]{%
  \expandafter\ifx\csname @seccntformat@#1\endcsname\relax
    \expandafter\@@seccntformat
  \else
    \expandafter
      \csname @seccntformat@#1\expandafter\endcsname
  \fi
    {#1}%
}
\newcommand*{\@seccntformat@section}[1]{%
  \S\csname the#1\endcsname\quad
}
\makeatother

\makeatletter
\renewcommand{\tocsection}[3]{%
  \indentlabel{\@ifnotempty{#2}{\S~\ignorespaces#1 #2.\quad}}#3}
\makeatother

\addto\captionsfrench{

}

\makeatletter
\@addtoreset{section}{part}% ou \numberwithin{section}{part}
\@addtoreset{proposition}{part}% ou \numberwithin{proposition}{part}
\@addtoreset{theoreme}{part}% ou \numberwithin{theoreme}{part}
\makeatother

\input cyracc.def 
 
\newtheorem{theoreme}{Théorème}
\newtheorem{proposition}{ Proposition}
\newtheorem{lemme}{ Lemme}
\newtheorem{corollaire}{Corollaire}

\newtheorem{definition}{Définition}

\theoremstyle{remark}

\def \1{\mathbb {1}}
\def \RM{\mathbb {R}}%        corps des reels
\def \NM{\mathbb{N}}%        entiers naturels
\def \ZM{\mathbb{Z}}%        entiers relatifs
\def \CM{\mathbb{C}}%        nombres complexes

\def \QM{\mathbb{Q}}%        nombres rationnels

 \def \Vol {{\rm Vol}}

\def \intr {{\rm int\,}}

\def \ord {{\rm ord\,}}

\def \Gr {{\rm Gr\,}}
\def \Der {{\rm Der\,}}

\def \Aut {{\rm Aut\,}}

\def \p {{\rm exp\,}}
\def \Id {{\rm Id\,}}

\def \d{\partial}%derivee partielle
\def\dt{\delta} 
\def\a{\alpha}
\def\b{\beta}
\def\e{\varepsilon}  
\def\g{\gamma}

\def\l{\lambda}

\def\p{\varphi}  
\def\lb{\left\{}
\def\rb{\right\}}
\def\G{\Gamma}   
\def\D{\Delta}
\def \s{\sigma}
\def \t{\tilde}

\def \to{\longrightarrow} 
\def \w{\wedge}
\def \alg{\mathfrak{g}}

\def\del{\nabla}
\def \< {{\langle }}
\def \> {{\rangle }}
\def \( {\left( }
\def \) {\right) }

\newcommand{\Bt}{{\mathcal B}}

\newcommand{\Dt}{{\mathcal D}}

\newcommand{\Lt}{{\mathcal L}}
\newcommand{\Mt}{{\mathcal M}}

\newcommand{\Ot}{{\mathcal O}}

\renewcommand{\mod}{{\rm  mod\,}}
\title[LA CONJECTURE DE HERMAN]{LA CONJECTURE DE HERMAN}
\author{  Mauricio  Garay}

\address{Max Planck,\ Institut für Mathematik\\
Vivatsgasse 7\\
53111 Bonn, Allemagne.}
\email{garay@mpim-bonn.mpg.de}
\begin{document}
\maketitle
%\noindent {\em  Les séries ne pourraient-elles pas, par exemple, converger quand $x^0_1$ et $x^0_2$ ont été choisies de telle sorte que le rapport
%$n_1/n_2$ soit incommensurable, et que son carré soit au contraire commensurable (\dots)?}
%\begin{flushright} H. Poincaré, {\em Méthodes mathématiques de la mécanique céleste.} \end{flushright}
%\bigskip
%\begin{flushright}{\em Cet article est dédié à Duco van Straten.}\end{flushright}
%\bigskip

\deactivateaddcontentsline
\section*{Introduction}

 La plupart des séries perturbatives de la physique mathématique divergent. La découverte -- due à Poincaré -- de ce phénomène pour les séries de Delaunay--Lindstedt en est probablement l'exemple le plus célèbre.  Poursuivant les travaux de Poincaré, Siegel montra que l'ensemble des fonctions analytiques pour lesquelles ces séries sont convergentes en un point critique est maigre au sens de Baire~\cite{Poincare_trois,Siegel_divergence,Siegel_Hamilton} (voir également~\cite{Perez_Marco}). 
 
 Parallèlement à ces résultats négatifs,  Poincaré et Siegel démon\-trèrent les premiers résultats de convergence sur les champs de vecteurs analytiques~\cite{Poincare_these,Siegel_vecteurs}.  Ces études aboutirent au théorème des tores invariants de Kolmogorov et à la théorie KAM~\cite{Arnold_KAM,Kolmogorov_KAM,Moser_KAM}. 
 
 Dans les années 90, Herman postula l'existence de tores invariants dans des systèmes hamiltoniens plus gé\-né\-raux que ceux obtenus par perturbation de systèmes intégrables. En 1998, lors du congrès international de Berlin, il formula la :\\
  
{\bf Conjecture.}\ \cite{Herman_ICM}  {\em   Un symplectomorphisme réel analytique  de fréquence diophantienne possède, au voisinage d'un point fixe elliptique, un ensemble de mesure positive de tores invariants.}\\

Herman conjectura également les deux variantes suivantes plus proches des résultats de Poincaré et de Siegel~:\\

{\em   Un  hamiltonien réel analytique  de fréquence diophantienne possède, au voisinage d'un  tore KAM, un ensemble de mesure positive de tores invariants.}\\

{\em   Un  hamiltonien réel analytique de fréquence diophantienne possède, au voisinage d'un point critique elliptique, un ensemble de mesure positive de tores invariants.}\\

Cet article a pour but de démontrer ces conjectures en contrebalançant les résultats de divergence de Poincaré et Siegel par des énoncés de type KAM. En fait, les trois dé\-monstrations ne diffèrent que par des détails de notations, je ne donnerai donc la dé\-mons\-tration que pour la troisième.

Le résultat de cet article est en fait plus fort que celui conjecturé par Herman~: sous des conditions dites de Bruno -- ce qui inclut le cas diophantien -- l'ensemble des tores invariants est, dans chacun des cas, paramétré un ensemble de densité égale à un en chacun de ses points.  
 
La démonstration de la conjecture de Herman est basée sur le principe suivant.  Dans un premier temps, nous montrerons que, sous des conditions de Bruno, un hamiltonien avec un point critique possède une variété lagrangienne complexe invariante par le flot hamiltonien.  Puis, dans un deuxième temps, nous verrons que l'existence d'une telle variété entraîne l'existence d'une famille de variétés lagrangiennes invariantes qui se concentrent au  point critique.  

Cette famille définit une déformation de la variété initiale. Mais alors que les déformations de la géométrie analytique ont gé\-né\-ralement pour base des espaces analytiques, celle-ci a pour base un ensemble totalement discontinu. Si l'on analyse la partie réelle de cette famille dans le cas elliptique, on voit que la  fibre spéciale est réduite à un point alors que celle des fibres gé\-né\-riques sont des tores. 

Ainsi, l'existence d'une variété lagrangienne invariante entraîne celle d'une famille de telles variétés. On pourrait se demander plus gé\-né\-ralement si ce principe se gé\-né\-ralise, mais cela semble hors de portée pour le moment.
 
  \noindent {\bf Remerciements.} { Ils vont tout d'abord à  H. Eliasson et à J. Féjoz qui m'ont aidé a concrétiser des idées demeurées pendant des années à l'état de conjectures.  Merci aussi à D. Kleinbock, J.-C. Yoccoz et B. Weiss  pour leurs explications sur l'approximation diophantienne, à R. Krikorian  et  B. Fayad pour leurs éclaircissements concernant la conjecture de Herman, à F. Jamet et R. Uribe pour leurs remarques concernant quelques aspects élémentaires de la théorie de la mesure et enfin à F. Aicardi pour le dessin qui illustre le théorème de densité arithmétique.
  
 Ce travail a été financé par le Max Planck Institut für Mathematik de Bonn et par le projet du Deutsche Forschungsgemeinschaft, SFB-TR 45, M086, {\em Lagrangian geometry of integrable systems.}}
\newpage 
   \tableofcontents
\reactivateaddcontentsline
%%%%%%%%%%%%%%%%%%%%%%%%%%%
%%%%%%%%%%%%%%%%%%%%%%%%%%%%%%
\section{\'Enoncé des résultats}
\subsection{Un résultat de densité sur les classes arithmétiques}
\label{SS::arithmetique}
  Notons $(\cdot,\cdot)$ le produit scalaire euclidien dans $\RM^n$.  
     Pour tout vecteur $\a \in \RM^n $, la suite numérique $\s(\a)$ définie par
$$\s(\a)_k :=\min \{ |(\a,i)|: i \in \ZM^n \setminus \{ 0 \}, \| i \| \leq 2^k \} $$
mesure l'éloignement du vecteur $\a$ au réseau des entiers. 
\begin{definition}  On appelle classe arith\-mé\-tique de $\RM^n$ associée à une suite $a=(a_k)$ l'ensemble
$$\Dt_a:=\{\a \in \RM^n: \s(\a)_k  \geq  a_k\} .$$
\end{definition}
 
 Soit $U \subset \RM^d$ un ouvert.
  \begin{definition}[\cite{Kleinbock,Kleinbock_Margulis}]\label{D::KM} Une application 
 $$g:U \to \RM^n$$
 de classe $C^k$ est dite  non-dégé\-né\-rée en un point $x \in \RM^d$ s'il existe un sous-espace vectoriel $P \subset \RM^n$ et un entier $l \leq k$ tel que~:
\begin{enumerate}[{\rm i)}]
 \item l'image de $g$ est contenue dans $g(x)+P$ ;
   \item les dérivées de $g$ d'ordre $l$ évaluées en $x$ engendrent $P$.
   \end{enumerate}
   \end{definition} 

  Pour $\a \in \RM^n$, nous noterons $B(\a,r)$ la boule de rayon $r$  centré $\a$, sans préciser la dimension de l'espace ambiant. Rappelons que la densité d'un sous ensemble mesurable $K \subset \RM^n$ en un point $\a$ est la limite (lorsqu'elle existe)~:
$$\lim_{r \to 0} \frac{\Vol(K \cap B(\a,r))}{\Vol(B(\a,r))}.$$

Si $u=(u_k)$ et $v=(v_k)$ sont deux suites numériques, on note $uv$ leur produit $(uv)_k=u_kv_k.$

\begin{theoreme} 
\label{T::arithmetique}
 Soit $a=(a_i),\rho=(\rho_i), \rho_i<1$ deux suites numériques décroissantes strictement positives et 
 $$f:\RM^d \to \RM^n $$
 une application non-dégé\-né\-rée.  Si
  $$\sum_{k \geq 0} (2^{(k+1)n+1} \sqrt{\rho_k}) <+\infty $$
  alors la densité de l'ensemble $f^{-1}(\Dt_{\rho a})$ est égale à $1$ en tout point de $f^{-1}(\Dt_a)$. 
\end{theoreme} 
 \begin{definition} Une suite numérique positive $(p_n)$ est à décroissance modérée si elle vérifie la condition
 $$-\sum_{n \geq 0}{\frac{ \log p_n'}{2^n}}<+\infty,\ p_n'=\min(1,p_n).  $$
 \end{definition}
De telles suites forment une algèbre~: le produit et la somme de deux suites à décroissance modérée est également à décroissance modérée.
 On dit que le vecteur $\a \in \RM^n$ {\em vérifie la condition de Bruno}  lorsque la suite $\s(\a)$ est à décroissance modérée~\cite{Brjuno}.

Une suite numérique positive $(p_n)$ est à {\em croissance modérée} si la suite $(p_n^{-1})$ est à décroissance modérée.
%%%%%%%%%%%%%%
\subsection{Forme normale d'un hamiltonien en un point critique de Morse}
Munissons $\RM^{2n}$ de coordonnés $q_i,p_i$ pour $i=1,\dots,n$ ainsi que de la forme symplectique standard
$$\omega:=\sum_{i=1}^n dq_i \w dp_i. $$
Désignons par $X,X'$ des sous-ensembles fermés de $\RM^n$. {\em Un symplectomorphisme $C^k$}
$$\p:X \to X' ,\ \p(X)=X'$$
est la restriction à $X$ d'un difféomorphisme de classe $C^k$ tel que
$$(\p^* \omega)_x=\omega,\ \forall x \in X'. $$
Plus généralement, un morphisme de Poisson d'un sous-ensemble d'une variété de Poisson  est la restriction d'un difféomorphisme de la variété qui préserve la structure de Poisson sur l'ensemble donné.

  Si $f,g$ sont des fonctions analytiques dans les variables $q,p$. Nous écrirons
$$f=g+o(l) $$
 lorsque la série de Taylor de $f-g$ est une somme de polynômes homogènes dont le degré est supérieur à $l$.

Si $H:(\RM^{2n},0) \to (\RM,0) $
est une fonction analytique avec un point critique de Morse en l'origine dont la partie quadratique est définie positive. 

Par un changement linéaire de variables qui préserve la forme symplectique, la fonction $H$ se met sous la forme
$$\sum_{i=1}^n \a_i (p_i^2+q_i^2)+o(2).$$ On dit alors que $H$ est {\em elliptique de fréquence}
 $$\a=(\a_1,\dots,\a_n) \in \RM^n.$$

Con\-si\-dé\-rons l'application
$$\pi:\RM^{2n} \to \RM^n,\ (q,p) \mapsto (p_1^2+q_1^2,p_2^2+q_2^2,\dots,p_n^2+q_n^2) $$
et notons $x_1,\dots,x_n$ les coordonnées sur $\RM^n$.
\begin{theoreme} 
\label{T::Herman}
Soit $b=(b_i)$ une suite numérique à décroissance modérée et 
$$H:(\RM^{2n},0) \to (\RM,0)$$ 
une fonction analytique avec un point critique elliptique de fréquence
$$\a=(\a_1,\dots,\a_n) \in \Dt_b .$$ 
Pour tout $k \geq 0$, il existe une application   $A:(\RM^n,0) \to \RM$ de classe $C^k$ et un germe de symplectomorphisme
de classe $C^k$~:
$$\p:(X,0) \to (X',0),\ X:=\pi^{-1} \circ (\del A)^{-1}(\Dt_b)$$
tels que
\begin{enumerate}[{\rm i)}]
\item $H \circ \p=A(p_1^2+q_1^2,\dots,p_n^2+q_n^2) $ ;
\item l'application gradient $\del A=(\d_{x_1}A,\d_{x_2} A,\dots,\d_{x_n} A):(\RM^n,0) \to (\RM^n,0)$ est non-dégé\-né\-rée ;
\item la restriction de $\p$ aux fibres de $\pi$ est analytique.
\end{enumerate}
\end{theoreme}
  Pour montrer la conjecture de Herman, on applique (par exemple) le deuxième théorème à la suite
$$b=\rho a,\ a:=\s(\a),\ \rho_i:=2^{-2(k+2)i+2}. $$
D'après le premier théorème, l'ensemble 
$$K:=\del A^{-1}(\Dt_b)$$
est de densité égale à un au point $\a \in \Dt_a$. L'ensemble $X$ est fibré au-dessus d'un ouvert dense de $K$ par des tores, il est donc également de densité égale à un en l'origine, de même que $X'$ qui est l'image de $X$ par un difféomorphisme.

%%%%%%%%%%%%%%%%
\subsection{Réseaux et classes arithmétiques}
 Au vecteur $\a \in \RM^n$, on associe le réseau $\G[\a]$ de $\RM^{n+1}$ définit par
$$\G[\a]:=\{(i,(\a,i)) \in \RM^{n+1}: i\in \ZM^n \} $$
où $(\cdot,\cdot)$ désigne le produit scalaire euclidien. 

Désignons par 
$$g_t:\RM^{n+1} \to \RM^{n+1} $$ 
l'application linéaire dont la matrice dans la base canonique est diagonale de coefficients~:
$$(e^{-t},e^{-t},\dots,e^{-t},e^t). $$

Soit $\G \subset \RM^{n+1}$ un réseau de l'espace euclidien de rang $n$. On pose
$$\dt(\G)=\inf_{\g \in \G} \| \g \| $$
où $\| \cdot \|$ désigne la norme euclidienne.

\begin{lemme}[cf. \cite{Dani,Kleinbock}]
\label{L::arithmetique}
Si $|(\a,i)| \leq a $ alors
$$\dt(g_t \G[\a]) \leq \e $$ où $\e,t$ sont définis par
$$\left\{ \begin{matrix}\e&=&\sqrt{2}\, e^t\, a\ ;\\ \e &=&\sqrt{2}\,e^{-t}\,\| i \|  . \end{matrix} \right.$$
\end{lemme}
\begin{proof}
Pour tout $x \in \RM^n$ et tout $y \in \RM$, on a~:
$$\| (x,y) \| \leq \sqrt{2} \max\left( \| x \| ,| y | \right) .$$
Par conséquent l'inégalité $|(\a,i)| \leq a $ entraîne~:
$$| g_t (i,(\a,i) | \leq  \sqrt{2} \max\left(e^{-t} \| i \| , e^t a \right)=\e $$
\end{proof}
La résolution explicite des équations du lemme donne
$$\left\{ \begin{matrix}\e&=&\sqrt{2a\| i \|}\, ; \\ t &=&\frac{1}{2}\log \frac{\| i \|}{a}   \end{matrix} \right. $$   
 
 \begin{theoreme}[\cite{Kleinbock_Margulis}, Proposition 3.4, Theorem 5.4]
 \label{T::KM} Soit $U \subset \RM^d$ un voisinage de l'origine et
$$f: U \to \RM^n$$
une application non-dégé\-né\-rée en l'origine. Il existe une constante $C>0$ telle que pour toute boule $B(0,r) \subset \RM^d$ de rayon suffisamment petit et pour tout $t>0$, on ait~:
$$\Vol\left( \{ x \in B(0,r): \dt\left(g_t \G[f(x)]\right)  \leq \e \} \right) \leq C\e^{1/dl} \Vol(B) $$
où $l$ est l'ordre de dérivées nécessaires pour engendrer l'image de $f$~(voir Définition~\ref{D::KM}).
\end{theoreme}
L'énoncé original de \cite[Proposition 3.4]{Kleinbock_Margulis}   ne concerne que les applications non-dégé\-né\-rées dont l'image n'est contenue dans aucun hyperplan. Dans notre situation, on commence par appliquer la proposition à la fonction 
$$\t f:U \to f(0)+P, x \mapsto f(x),\ P=f(U).$$
La conclusion de la proposition étant valable pour $\t f$, elle l'est également pour $f$~(voir aussi~\cite{Kleinbock}).
%%%%%%%%%%%%%%%%%%%%%%%%%%%%%%%%%%%%%%%%%%%%%%
\subsection{Démonstration du Théorème \ref{T::arithmetique}}
Notons $[\cdot ]$ la partie entière et considérons l'application
$$\p:\ZM^n \to \NM,\ i \mapsto \left[{\log}_2 \| i \| \right]+1. $$  
Pour $i \in \ZM^n$, $\p(i)$ est le plus petit entier tel que $i$ soit contenu dans la boule  de rayon $2^{\p(i)}$ centré en l'origine.

Fixons $i \in \ZM^n $ et posons $k:=\p(i)$. L'ensemble
$$ M_i:=\{\b \in  \RM^{n}: \left| (\b,i) \right| <  \rho_k a_k \} $$
 est  une bande de largeur $\rho_k a_k/\| i \|$ et la réunion des $M_i$ prise sur tous les $i \in \ZM^n$ est le complémentaire de la classe arith\-mé\-tique $\Dt_{\rho a}$.
 
\begin{figure}[ht]
\centerline{\epsfig{figure=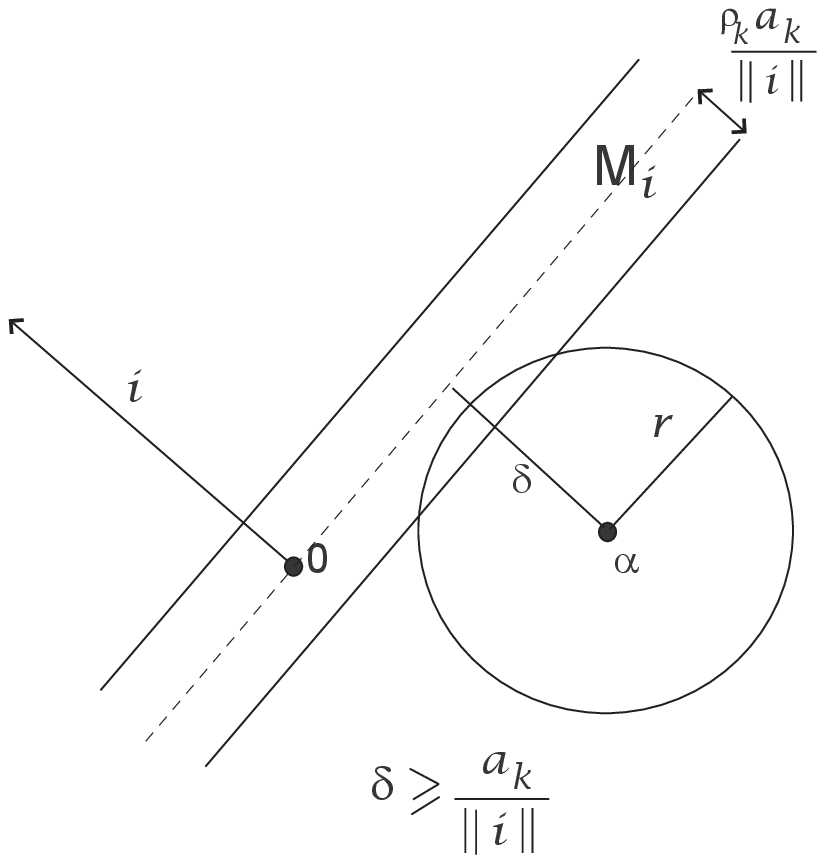,height=0.52\linewidth,width=0.48\linewidth}}
\end{figure} 
 Soit $\a$ un vecteur de  $ \Dt_a$. Notons $\dt$ la distance de $\a$ à l'hyperplan orthogonal au vecteur $i$.
 L'ensemble $M_i$  ne peut intersecter la boule $B(\a,r)$  que si
 $$r > \dt-\frac{\rho_k a_k}{\| i \|}.  $$
 Comme $\a \in \Dt_a$, on a donc nécéssairement
 $$\dt > \frac{a_k}{\| i \|} \geq \frac{a_k}{2^k}  .$$
Ce qui donne finalement
$$ \frac{(1-\rho_k)a_k}{2^k} < r.$$

Comme la suite $\rho$ est sommable, il existe $N$ tel que
$$\rho_k<\frac{1}{2},\ \forall k \geq N. $$
Choisissons 
$$r <\inf \lb \frac{(1-\rho_k)a_k}{2^k}: k \leq N \rb.$$ On a alors 
$$ \p(i)<N \implies M_i \cap B(\a,r)=\emptyset.$$ 

Ceci montre que si $M_i$ intersecte la boule $B(\a,r)$ alors le vecteur $i \in \ZM^n$ doit être un élément  de l'ensemble
 $$I_r:=\{ i \in \ZM^n: \frac{a_k}{2^{k+1}} < r,\ k =\p(i) \}. $$
 
 Posons
$$\left\{ \begin{matrix} &\e_k &:=&\sqrt{2^{k+1} a_k \rho_k}\ ;\\ \\ & t_k&:=&\displaystyle{\frac{1}{2}\log \frac{\| i \|}{a_k}}. \end{matrix} \right.$$ 
   D'après le lemme \ref{L::arithmetique}, on a~:
 $$f(x) \in M_i \implies  \dt(g_{t_k}\G[f(x)]) \leq  \e_k . $$
Par conséquent, le théorème \ref{T::KM} entraîne l'existence de constante $C,\g>0$ telles que
 $$\Vol (B(0,r) \cap f^{-1} (M_i)) \leq C\e_k^\g Vol(B(0,r)). $$
  
 L'application
 $$f:\RM^d\to \RM^n $$
 étant différentiable, d'après le théorème des accroissements finis, il existe une constante $\kappa$ telle que
 $$f(B(0,r)) \subset B(\a,{\kappa r}),\ f(0)=\a, $$
 pour tout $r$ suffisamment petit. Pour tout $r$ suffisament petit, on a donc~:
 $$f^{-1} (M_i) \cap B(0,r) \neq \emptyset \implies i \in I_{\kappa r}.$$
 
Ceci montre que  le complémentaire de $f^{-1}(\Dt_{\rho a})$ dans la boule $B(0,r)$ a sa mesure majorée par
 $$C\, \Vol(B(0,r))\left(\sum_{i \in I_{\kappa r}}\sqrt{2^{\p(i)+1} a_{\p(i)} \rho_{\p(i)}}\right)^\g.$$
 
 Par définition de $I_{\kappa r}$, on a~:
 $$\sum_{i \in I_{\kappa r}}\sqrt{2^{{\p(i)}+1} a_{\p(i)} \rho_{\p(i)}} < 2\sqrt{ r}\sum_{i \in I_{\kappa r}} 2^{\p(i)}\sqrt{\rho_{\p(i)}} .$$
 et
$$\sum_{i \in \ZM^n} 2^{\p(i)}\sqrt{\rho_{\p(i)}} =\sum_{k \geq 0} \sum_{\p(i)=k} 2^{\p(i)}\sqrt{\rho_{\p(i)}} . $$
Notons $\#-$ le cardinal. On a~:
$$\#\{ \p(i)=k \}=\#\{ \p(i)\leq k\}-\#\{ \p(i)\leq k-1\} \leq 2^{(k+1)n} $$
donc 
$$\sum_{i \in \ZM^n} 2^{\p(i)}\sqrt{\rho_{\p(i)}} \leq \sum_{k \geq 0} 2^{(k+1)n+k}\sqrt{\rho_k} .$$

Par hypothèse, la série du membre de droite est convergente. Ceci montre que les sommes
$$ \sum_{i \in I_{\kappa r}}\sqrt{2^{\p(i)+1} a_{\p(i)} \rho_{\p(i)}}$$  
tendent vers $0$ avec $r$. Le théorème est démontré. 

%%%%%%%%%%%%%%%%%%%%%%%%%%%%%%%%%%%%%%
  %%%%%%%%%%%%%%
%%%%%%%%%%%%%%%%%%%%%%%%%%%%

 \section{Espaces vectoriels échelonnés}
  %%%%%%%%%%%%%%%%%%%%%%%%%%%%%%
  %%%%%%%%%%%%
 %%%%%%%%%%%%%%%%%%%%%%%%%%%%%
\subsection{Définition}
\label{SS::definition}
Une {\em $S$-échelle de Banach} est une famille décroissante d'espaces de Banach $(E_s)$, $s \in ]0,S[$, telle que 
les inclusions 
$$E_{s+\s} \subset E_s,\ s \in ]0,S[,\ \s \in ]0,S-s[$$ 
soient de norme au plus~$ 1$.

Soit $E$ un espace vectoriel topologique. Un {\em $S$-échelonnement} de $E$ est une échelle  $(E_s)$
de sous-espaces de Banach de $E$ telle que
\begin{enumerate}[{\rm i)}] 
\item  $E=\bigcup_{s \in ]0,S[} E_s$ ;
\item la topologie  limite directe induite par les inclusions $E_s \subset E$ coïncide avec celle de $E$.
\end{enumerate}
(Rappelons que si $f_s:X_s \to X$ un famille d'application d'espaces topologiques $(X_s)$ dans un ensemble $X$. On appelle topologie {\em  limite directe} sur $X$, la topologie  la plus fine sur $X$ qui rend les applications $f_s$ continues~:
 $$U \subset X {\rm \ est\ ouvert\ }\ \iff  f^{-1}_s(U) {\rm \  est\ ouvert\  dans}\ X_s,\ {\rm pour\ tout}\ s.)$$
 
 L'intervalle $]0,S[$ s'appelle {\em l'intervalle d'échelonnement}.  Si $F$ est un sous-espace vectoriel fermé d'un espace vectoriel échelonné $E$ alors $E/F$ est échelonné par les espaces de Banach
$$(E/F)_s:=E_s/(E \cap F)_s.$$

Finalement, rappelons que le produit tensoriel topologique de deux espaces de Banach $E,F$, noté $E \hat \otimes F$ est le complété de $E \otimes F$ pour la norme
$$\| \sum_{i=1} x_i \otimes y_i \| =\inf \sum_{i=1} \|  x_i \|  \otimes \|  y_i \| $$
où la borne inférieure est prise sur les différentes écritures possibles du produit tensoriel~\cite{Schatten} (voir aussi \cite{Grothendieck_PTT}).  

On définit ainsi le produit tensoriel de deux échelonnements
$$(E \hat \otimes F)_{s}:=E_s \hat \otimes F_s. $$
La limite directe des espaces de Banach $(E \hat \otimes F)_{s}$ définit un espace vectoriel topologique que nous noterons $E \hat \otimes F$. Dans les cas que nous allons considérer, cet espace ne dépend pas du choix des échelonnements et notre définition est équivalente à celle de Grothendieck~\cite{Grothendieck_PTT}. 
 
 % %%%%%%%%%%%%%%%%%%%%%%%%%%%%%%%%%%%%%%%%%
\subsection{Morphismes d'un espace vectoriel échelonné}
\label{SS::morphismes}
Soit $E,F$ deux espaces vectoriels $S$-échelonnés.

Nous dirons d'une application linéaire que c'est  un {\em morphisme} entre des espaces vectoriels échelonnés $E,F$, si pour tout $s \in ]0,S[$, 
il existe $s' \in ]0,S[$ tel que l'espace de Banach $E_{s}$ est envoyé continûment dans $F_{s'}$. Nous avons ainsi définit la  {\em catégorie des espaces vectoriels échelonnés}.

Nous désignerons par $\Lt(E,F)$ l'espace vectoriel des morphismes de $E$ dans $F$ et lorsque $E=F$, nous utiliserons la notation $\Lt(E)$ au lieu de $\Lt(E,E)$. Il n'y pas de raison, a priori, pour que $\Lt(E,F)$ coïncide avec l'espace des applications linéaires continues de $E$ dans $F$, mais dans les exemples concrets que nous allons traiter ce sera toujours le cas.
 
   Si $\| \cdot \|$ désigne la norme d'opérateur sur l'espace de Banach $\Lt(E_{s'},F_s)$, nous noterons $\| u \|$ la norme de l'opérateur défini par restriction de $u$ à $E_{s'}$.
Venons-en à la notion de convergence d'une suite de morphismes. La norme d'opérateur induit sur  les espaces vectoriels $\Lt(E_{s'},F_s)$, une structure d'espace de Banach. 
\begin{definition}Une suite de morphismes $(u_n)$ de $\Lt(E,F)$ converge vers un morphisme $u \in \Lt(E,F)$ si pour tout $s' \in ]0,S[$, il existe $s \in ]0,S[$ tel que la restriction de $(u_n)$ définisse une suite de $\Lt(E_{s'},F_s)$ qui converge vers la restriction de $u$.
\end{definition}
Un sous-ensemble $X$ de $\Lt(E,F)$ sera dit {\em fermé} si toute suite convergente de points de $X$ à sa limite dans $X$. (L'utilisation du mot «fermé» est légèrement abusive, car il ne s'agit pas a priori du complémentaire d'un ouvert.)
%%%%%%%%%%%

%%%%%%%%%%%%%%%%%%%%%%%%%%%%%%%%%%%%
\subsection{Filtration d'un espace vectoriel échelonné}
Soit $E$ un espace vectoriel échelonné.  Les sous-espaces vectoriels
 $$E^{(k)}=\{ x \in E: \exists C,\tau,\ | x|_s \leq Cs^k,\ \forall s \leq \tau \},\ k \geq 0$$
 filtrent l'espace vectoriel $E$~:
 $$E:=E^{(0)} \supset E^{(1)}\supset E^{(2)} \supset \cdots $$
 Par ailleurs, on pose $E^{(-k)}=E$ pour tout $k \geq 0$.
 
\begin{definition} On appelle ordre de $x \in E$, noté $ord(x)$ le plus grand $j$ tel que
 $x \in E^{(j)}.$
\end{definition}  

Le {\em gradué, noté $\Gr(E)$, associée à un espace vectoriel échelonné $E$} est l'espace vectoriel échelonné
$$\Gr(E):=\bigoplus_{i \geq 0} E^{(i)}/E^{(i+1)}.  $$
Si $M$ est un sous-ensemble de $E$, on désigne par $\Gr(M)$ le plus petit sous-espace vectoriel de $\Gr(E)$ qui contient l'image de $M$
par la projection canonique
$$E \to \Gr(E). $$

 %%%%%%%%%%%%%%%
 \subsection{  $\tau$-morphismes, morphismes bornés}

Conservons les notations du chapitre précédent. 
\begin{definition} Un morphisme $ u \in \Lt(E,F)$  est appelé  un $\tau$-morphisme, $\tau<S$,
si pour tout $s'  \in ]0,\tau]$ et pour tout $s \in ]0,s'[$, on a l'inclusion $ u(E_{s'}) \subset F_s$ et $ u$ induit par restriction une application linéaire continue
$$u_{s',s}~:~E_{s'}~\to~F_s.$$
\end{definition}
On a alors des diagrammes commutatifs
$$\xymatrix{ & \ F_s \ar@{^{(}->}[d]\\
 E_{s'} \ar[r]^-{u_{\mid E_{s'}}} \ar[ru]^{u_{s',s}}  & F }
$$
pour tout $s'  \in ]0,\tau]$ et pour tout $s \in ]0,s'[$, la flèche verticale étant donnée par l'inclusion $F_s \subset F$.
\begin{definition}[\cite{groupes}]
\label{D::borne}
  Un $\tau$-morphisme $ u:E \to F$ d'espaces vectoriels $S$-échelonnés est dit 
    $k$-borné, $k \geq 0 $ s'il existe un réel $C>0$ tel que~:
  $$| u(x) |_s \leq C \sigma^{-k} | x |_{s+\sigma},\ {\rm pour\ tous\ }\ s \in ]0,\tau[,\ \s \in ]0,\tau-s],\ x\in E_{s+\s} . $$
\end{definition}
Un morphisme est dit {\em $k$-borné} (resp. {\em borné}) s'il existe $\tau$ (resp. $\tau$ et $k$) pour lequel (resp. lesquels) c'est un $\tau$-morphisme $k$-borné. Lorsque $E=E_s$ et $F=F_s$ sont des espaces de Banach, on retrouve la définition habituelle de morphismes bornés.
(Nous n'utiliserons pas la notion plus gé\-né\-rale d'application linéaire bornée d'un espace localement convexe, notre terminologie ne devrait donc pas porter à confusion.)

Si  $u:E \to F$ est $k$-borné alors $u(E^{(j+k)}) \subset F^{(j)} .$

L'espace vectoriel des $\tau$-morphismes (resp. des morphismes)  $k$-bornés
entre $E$ et $F$ sera noté $\Bt^k_\tau(E,F)$ (resp. $\Bt^k(E,F)$). On note $N_\tau^k(u)$ la plus petite constante~$C$ vérifiant l'inégalité de la définition \ref{D::borne}.

 On vérifie facilement que  si $E,F$ sont des espaces vectoriels $S$-éche\-lon\-nés alors les espaces vectoriels normés
$(\Bt^k_\tau(E,F),N_\tau^k),\ \tau \in ]0,S[$, forment une $S$-échelonnement de $\Bt^k(E,F)$.  

La propriété pour un endomorphisme {\em surjectif} d'être borné passe au quotient~: tout endomorphisme $k$-borné surjectif $ u:E \to E$ définit un endomorphisme $k$-borné sur l'espace quotient $E/F$, pour tout sous-espace vectoriel fermé $F \subset E$. 

\begin{definition} Un sous espace vectoriel fermé $F$ d'un espace vectoriel échelonné $E$ est dit  $m$-direct (ou tout simplement  direct) si pour tout $s$ et pour tout $n \geq 0$, il existe un supplémentaire de $F_s^{(n)}$ dans $E^{(n)}_s$ et si la projection sur $F$ est $m$-bornée de norme $1$. 
\end{definition}

 %%%%%%%%%%%%%%%%%%%
 \subsection{\'Echelonnement des germes de fonctions holomorphes}
 \label{SS::echelon}
 Munissons  l'espace vectoriel des germes de fonctions holomorphes à l'origine  $E:=\Ot_{\CM^n,0}$ de la topologie de la convergence uniforme sur les compacts de $\CM^n$. 
 
 Fixons $S>0$ et $s \in ]0,S[$. Notons $E_s$ l'espace vectoriel des fonctions continues sur le polydisque 
 $$D_s:= \{ z \in \CM^i: \sup_{i=1,\dots,n }| z_i | \leq s \}$$
 qui sont holomorphes dans l'intérieur de ce polydisque. La norme
 $$| f |_s:=\sup_{z \in D_s} |f(z)| $$
muni l'espace vectoriel $E_s$ d'une structure d'espace de Banach et les $(E_s)$ forment un échelonnement de $E$~\cite{Douady_these,Nagumo}.  

Notons $\G(-,-)$ le foncteur des sections globales et  $\intr(-)$ l'intérieur, on a donc~:
$$E_s:= \G(\intr(D_s),\Ot_{\CM^n}) \cap C^0(D_s,\CM).$$
Les inégalités de Cauchy montrent que pour cet échelonnement,  tout opérateur différentiel d'ordre $k$ sur $\Ot_{\CM^n,0}$ définit un morphisme $k$-borné.

Par unicité du complété, la multiplication
$$\Ot_{\CM^m,0} \hat \otimes \Ot_{\CM^n,0} \to \Ot_{\CM^{m+n},0}, \sum_{i \geq 0} f_i \otimes g_i \mapsto \sum_{i \geq 0} f_i g_i $$
est un isomorphisme d'espaces vectoriels topologiques et nous identifierons souvent les produits tensoriels avec leur image par multiplication, lorsque celle-ci est injective.

Notons $\Mt_{\CM^n,0}$ l'idéal maximal de $\Ot_{\CM^n,0}$~:
$$\Mt_{\CM^n,0}=\{ f \in \Ot_{\CM^n,0}:f(0)=0 \}. $$
On a 
$$(\Ot_{\CM^n,0})^{(k)}=\Mt_{\CM^n,0}^k .$$
Par conséquent, la filtration d'espace vectoriel échelonné coïncide avec celle donnée par les puissances de l'idéal maximal.

Le gradué associé à cette filtration est un espace vectoriel isomorphe à un espace de polynômes
$$ \Gr(\Ot_{\CM^n,0})=\bigoplus_k \Mt_{\CM^n,0}^k/\Mt_{\CM^n,0}^{k+1}\approx \CM[z_1,\dots,z_n]$$
que l'on peut identifier à un sous-espace vectoriel de $\Ot_{\CM^n,0}$.

Un autre échelonnement $(H_s)$ de $\Ot_{\CM^n,0}$ s'obtient en remplaçant l'espace des fonctions continues par les fonctions $L^2$~:
$$E_s:= \G(\intr(D_s),\Ot_{\CM^n}) \cap L^2(D_s,\CM).$$
Comme toute fonction continue sur un compact est intégrable, l'identité définit un morphisme $0$-borné
$$J :(E_s) \to (H_s) $$ 
entre l'espace vectoriel $E$ échelonné par les $(E_s)$ dans le même espace vectoriel échelonné par les $(H_s)$.
\begin{proposition} L'inverse de l'application $J$ est un morphisme $1$-borné.
\end{proposition}
Notons $| \cdot |_s$ la norme de $H_s$ et $dV$ le volume euclidien sur $\CM^n \approx \RM^{2n}$. 
Pour $z \in D_s$ et $\s \geq 0$ fixés, on pose
$$f(z+\s)=\sum_{j \geq 0} a_j \s^j,\ a_j \in \CM. $$
On a
$$|f(z) |^2 \s^2= |a_0|^2 \s^2  \leq \sum_{j \geq 0} |a_j|^2 \s^{2j+2}= \int_{z+D_\s } | f(z)|^2 dV . $$
Comme $z+D_\s \subset D_{s+\s}$, on a l'estimation
$$\int_{z+D_\s } | f(z)|^2 dV   \leq \int_{D_{s+\s} } | f(z)|^2 dV= | f |_{s+\s}^2,$$
d'où la proposition.  
 
Cette proposition implique que tout sous-espace vectoriel de $E=\Ot_{\CM^n,0}$ muni de l'échelonnement $(E_s)$ est $1$-direct. 
 Ces considérations s'étendent sans difficultés aux germes le long d'un compact de $\CM^n$.

 %%%%%%%%%%%%%%%%%%%%%%%%%%%%%
\subsection{L'exponentielle d'un morphisme borné}
\label{SS::produits}
 Si $u,v$ sont des morphismes,  respectivement $k$ et $k'$ borné, alors leur composition $u  v$ est
  $(k+k')$-borné et on a l'inégalité
  $$N_\tau^{k+k'}( u v) \leq 2^{k+k'} N_\tau^k(u)N_\tau^{k'}(v) .$$
  En effet~:
 $$| (u v)(x) |_s \leq N_\tau^k(u) \frac{2^k}{\sigma^{k}} |v (x)|_{s+\sigma/2} \leq   N_\tau^k(u)N_\tau^{k'}(v) \frac{2^{k+k'}}{\sigma^{k+k'}} |x|_{s+\sigma} $$
 pour  tout $x\in E_{s+\s}$. Plus gé\-né\-ralement, on a la
\begin{proposition}
\label{P::morphismes}
Le produit de $n$ morphismes $k_i$ bornés $u_i,\ i=1,\dots,n$, est un morphisme $k $-borné avec $k:=\sum_{i=1}^n k_i$ et plus précisément
  $$N_\tau^k( u_1  \cdots  u_n) \leq n^k \prod_{i=1}^n N_\tau^{k_i}(u_i),\ .$$
  De plus, si $u_1=\dots=u_n=u$ est d'ordre $1$ alors
   $$\frac{N_\tau^n( u^n )}{n!} \leq 3^n  N_\tau^1(u)^n. $$
 \end{proposition} 
 \begin{corollaire}
 Soit $u$ un $\tau$-morphisme $1$-borné. Si l'inégalité $3N_s^1(u) ~<~s $ est satisfaite pour tout $s \leq \tau$ alors
 la série 
 $$e^u:=\sum_{j \geq 0}\frac{u^j}{j!} $$
 converge vers un morphisme de $E$, et plus précisément
$$| e^u x |_{\l s} \leq \sum_{j \geq 0} \frac{(3N^1_s(u))^j }{(1-\l)^j s^j} | x |_s=\frac{1}{1-\frac{3N^1_s(u) }{(1-\l) s} }   |x |_s$$
pour tous $\l \in ]0, 1-\frac{3N^1_s(u)}s[$, $s \in ]0, \tau]$ et $x \in E_s$.  
\end{corollaire}
Nous dirons qu'un $\tau$-morphisme $u$ est {\em exponentiable} si pour tout $s \leq \tau$, le $\tau$-morphisme $u$ vérifie l'inégalité $3N_s^1(u) ~<~s $.

 Finalement, remarquons que deux morphismes 1-bornés $u,v \in \Bt^1(E)$ exponentiables qui commutent vérifient l'égalité
 $$e^{u+v}=e^{u}e^{v}.$$  

\begin{proposition} 
\label{P::produits}
Soit $E$ un espace vectoriel échelonné.   Soit $(u_n) \subset \Bt^1_\tau(E)$une suite de $\tau$-morphismes $1$-bornés exponentiables.
Si la série numérique
$$\sum_{n \geq 0} N_\tau^1(u_n) $$
est convergente alors la suite $(g_n)$ définie par
$$g_n:=e^{u_n}e^{u_{n-1}}\cdots e^{u_0}  $$
converge vers un élément inversible de $\Lt(E)$.
 \end{proposition}
 Pour démontrer ce résultat, commençons par gé\-né\-raliser le corollaire précédent. 
 \begin{lemme} 
\label{L::produits}
Soit $(u_n)$une suite de $\tau$-morphismes $1$-bornés exponentiables.
Pour tout $s \leq \tau$, la norme du morphisme
$$g_n:=e^{u_n}e^{u_{n-1}}\cdots e^{u_0}  $$
vérifie l'inégalité
$$| g_n x |_{\l s} \leq \left( \prod_{i=0}^n \frac{1}{1-\frac{3}{(1-\l)s}N_s^1(u_i)}\right) | x |_s $$
pourvu que $\l$ vérifie
$$\max_{i \leq n}\frac{3}{(1-\l)s}N_s^1(u_i)<1. $$
 \end{lemme}
\begin{proof}
Notons $\D_j \subset \ZM^j$ les suites $ i=(i_1,\dots,i_j)$ dont les éléments sont dans l'ensemble
 $\{0,\dots,n \}$  et telles que $i_p \geq i_{p+1}$.  On a alors la formule
 $$\prod_{i=0}^n \frac{1}{1-z_i}=\sum_{j \geq 0} \sum_{i \in \D_j}z^i,\ z^i:=z_1^{i_1} z_2^{i_2} \dots  z_j^{i_j}$$
 et plus gé\-né\-ralement
 $$\prod_{i=0}^n \frac{1}{1-\a z_i}=\sum_{j \geq 0} \sum_{i \in \D_j}\a^ j z^i.$$
 
On pose
$$u[i]:=u_{i_1} u_{i_2} \cdots u_{i_j},\ i \in \D_j. $$
Développons $g_n$ en série puis regroupons les termes suivant l'ordre en $t$, il vient~:
$$g_n=\sum_{j \geq 0} (\sum_{i \in \D_j} u[i]) \frac{t^j}{j!}=1+(\sum_{i=0}^n u_i)t+(\sum_{i=0}^n u_i^2+\sum_{j=0}^n \sum_{i=j+1}^n u_i u_j) \frac{t^2}{2}+\dots .$$
Posons 
$$z_{i,s}:=N_s^1(u_{i})\ {\rm et} \ N_s^j(u[i]):=N_{\tau}^j(u_{i_1} u_{i_2} \cdots u_{i_j}).$$
De la proposition \ref{P::morphismes}, on déduit l'inégalité~:
$$\frac{1}{j!}N_s^j( u[i])  \leq 3^j\prod_{p=0}^jN_s^1(u_{i_p})= 3^j z^i_s,\ i \in \D_j$$
et par suite
$$\left| u[i] (x)\right|_{\l s} \leq   \left(\frac{ 3}{(1-\l)s} \right)^j z^i_s  | x|_s,\ \forall \l \in ]0,1[. $$
On obtient bien
$$| g_n x |_{\l s} \leq  \left( \sum_{j \geq 0}\sum_{i \in \D_j} \a^j z^i_s \right)| x|_s,\ \a=\frac{3}{(1-\l)s} .  $$
Ce qui démontre le lemme.
 \end{proof}
 
 Achevons la démonstration de la proposition. Pour cela, fixons $s \in ]0,\tau]$ et choisissons $\l \in]0,1[$ tel que
 $$ \sup_{n \geq 0}\frac{3}{(1-\l)s}N_s^1(u_n)) <1.$$
 
 Il est possible de trouver de tels $\l$ car la suite  $(N_s^1(u_n))$ étant majorée par $(N_\tau^1(u_n))$, elle
 tend vers $0$ quand $n$ tend vers l'infini.

 Montrons tout d'abord que la suite $(g_n)$ définit, par restriction, une suite uniformément bornée d'opérateurs dans $\Lt(E_s,E_{\l s})$. 
Pour cela, notons $\| \cdot \|_{\l}$ la norme d'opérateur dans $\Lt(E_s,E_{\l s})$. Le lemme précédent donne l'estimation
$$ \| g_n \|_\l \leq   \prod_{i=0}^n \frac{1}{1-\frac{3}{(1-\l)s}N_s^1(u_i)}.$$

En prenant le logarithme du membre de droite, on voit que le produit converge quand $n$ tend vers l'infini, car la série de terme gé\-né\-ral
$(N_s^1(u_n)) $
est convergente. Comme chacun des facteurs de ce produit est au moins égal à un, on obtient l'inégalité
$$  \| g_n \|_\l \leq  C_\l,\ C_\l:= \prod_{i \geq 0} \frac{1}{1-\frac{3}{(1-\l)s}N_s^1(u_i)}.$$
Ce qui démontre l'assertion.

Soit à présent, $\mu \in ]0,1[$ vérifiant l'inégalité
 $$ \sup_{n \geq 0}\frac{3}{(1-\mu)\l s}N_{\l s}^1(u_n)) <1.$$ Nous allons montrer que la suite $(g_n)$ définit, par restriction, une suite de Cauchy dans $\Lt(E_s,E_{\mu \l s})$, la proposition en découlera. 

Je dis que la série de terme gé\-né\-ral $ \| g_n-g_{n-1}\|_{\l \mu }$ est convergente. Pour le voir, écrivons
$$g_n-g_{n-1}=(e^{u_n}-\Id)g_{n-1} $$
où $\Id \in \Lt(E)$ désigne l'application identité.

En développant l'exponentielle en série, on obtient l'inégalité~:
$$ | (e^{ u_n}-\Id) y |_{\l \mu s} \leq  \left( \sum_{j \geq 0} \frac{ ( 3 N_{\l s}^1(u_n))^{j+1} }{((1-\mu) \l s )^{j+1} } \right) | y |_{\l s}=
\frac{3}{1-\mu-\frac{3N^1_{\l s}(u_N)}{\l s}} \frac{N^1_{\l s}(u_n)}{\l s} | y |_{\l s} ,$$
pour tout $y \in E_{\l s}$. En prenant $y=g_n x$, ceci nous donne l'estimation
$$ \| (e^{ u_n}-\Id) g_{n-1} \|_{\l \mu} \leq   \frac{3 C_\l}{1-\mu-\frac{3N^1_{\l s}(u_n)}{\l s}} \frac{N^1_{\l s}(u_n)}{\l s}.$$
 
La quantité 
$$K_{\l,\mu}:=\sup_{n \geq 0} \frac{3 C_\l}{1-\mu-\frac{3N^1_{\l s}(u_n)}{\l s}} $$
est finie car la suite $3N^1_{\l s}(u_n)$ tend vers $0$ lorsque $n$ tend vers l'infini. 
Nous avons donc montré l'estimation
$$\| g_n-g_{n-1} \|_{\l \mu} \leq K_{\l,\mu} \frac{N^1_{\l s}(u_n)}{\l s}.$$
Il ne nous reste plus qu'à utiliser l'inégalité triangulaire pour voir que $(g_n)$ définit une suite de Cauchy de l'espace de Banach $\Lt(E_s,E_{\mu \l s})$~:
$$\| g_{n+p}-g_{n}  \|_{\l \mu } \leq \sum_{i = 1}^{p}  \| g_{n+i}-g_{n+i-1}\|_{\l \mu  } \leq K_{\l,\mu}\left(\sum_{i = 1}^{p} \frac{3N^1_{\l s}(u_{n+i})}{\l s}\right).$$
Nous avons donc montré que la suite $(g_n)$ converge vers un élément $g \in \Lt(E)$. On démontre de même que la suite $(h_n)$  définie par
$$h_n=e^{- u_0}e^{- u_1} \cdots e^{-u_n}$$ 
converge vers un élément $h \in \Lt(E)$. Pour tout $n \in \NM$, on a~:
$$g_nh_n=h_ng_n=\Id$$
 donc $gh=hg=\Id$. Ce qui montre que $h$ est l'inverse de $g$.
 La proposition est démontré.
        
%%%%%%%%%%%%%%%%%%%%%%%%%%%
 \subsection{Morphismes modérés}
 Soit $S>0$ et $E:=(E_n)$ $F:=(F_n)$ des suites d'espaces vectoriels $S$-échelonnés. Nous noterons abusivement de la même façon les normes
 pour les différentes valeurs de $n$. Un {\em morphisme} $u:E \to F$  est une collection de morphismes d'espaces vectoriels échelonnés.
 
Les notions définies pour les espace vectoriels échelonnés s'étendent naturellement aux suites d'espaces vectoriels échelonnés. Par exemple, si  les $u_n$ sont des $\tau$-morphismes (resp. des $\tau$-morphismes $k$-bornés) nous dirons
que $u$ est un $\tau$-morphisme (resp.  $\tau$-morphisme $k$-borné).  En accord avec notre convention,  nous notons $N^k_\tau$ la norme de l'espace vectoriel $\Bt^k(E_n,F_n)$ sans préciser l'indice $n$.

 \begin{definition}
 Un  $\tau$-morphisme
 $$u=(u_n):E \to F $$
 est dit $k$-modéré si, pour tout $s$, la suite numérique $N^k_{\tau}(u_n)$, appelée la norme de $u$, est à croissance modérée.
 \end{definition}
 Nous noterons $\Mt^k(E,F)$ l'ensemble des morphismes $k$-modérés.

\begin{definition} Un quasi-inverse à droite d'un morphisme
$$u:E \to F $$
est un morphisme
$$v:F \to E $$
tel que
$$u_n \circ v_n(x)=x\ \mod F^{(n+1)} $$
pour tout $n \geq 0$.
\end{definition}
On définit de manière similaire les quasi-inverses à gauche. On pourra comparer cette notion à celle introduite par Moser et Zehnder~\cite{Moser_pde,Zehnder_implicit}.
 
   %%%%%%%%%%%%%%% 
  %%%%%%%%%%%%%%%%%%%%%%%%%%% 

 \subsection{Approximations d'espaces vectoriels échelonnés}

 Soit  $E_\infty$ un espace vectoriel $S$-échelonné.
\begin{definition} Une approximation de $E_\infty$  consiste en la donnée d'une suite d'espace vectoriels échelonnés $(E_n)$ ainsi que de deux suites de morphismes  $0$-bornés , appelés morphismes de restrictions, entre espaces vectoriels~: 
$$ E_0 \stackrel{r_0}{\to} E_1 \stackrel{r_1}{\to} E_2 \stackrel{r_2}{\to}\dots $$
et
$$E_n \stackrel{s_n}{\to} E_\infty $$
tels que 
\begin{enumerate}[{\rm i)}]
\item $s_{n+1} r_n=s_n $ ;
\item la norme des $r_n,s_n$ est au plus égale à $1$.
\end{enumerate}
\end{definition}

Le produit tensoriel topologique de deux approximations est définit en prenant le produit tensoriel des espaces vectoriels échelonnés~:
 $$(E \hat \otimes F)_{n,s}:=E_{n,s} \hat \otimes F_{n,s}. $$

   %%%%%%%%%%%%%
 \subsection{Approximations tautologiques}
 Soit $E$ un espace vectoriel $S$-échelonné. Consi\-dé\-rons la suite d'espaces vectoriels échelonnés $(E_n)$ définie de la façon suivante. Comme espaces vectoriels topologiques, les $E_n$ sont tous identiques~: 
$$E_n=E {\rm\ et \ } E_\infty=E.$$ 
En revanche, l'échelonnement de $E_n$ diffère en fonction de $n$~:
$$(E_n)_s :=E_{s_n}, (E_\infty)_s=E_s$$
avec
$$ s_n:=\frac{n+2}{n+1}\,s,\  s \in \left]0,\frac{S}{2}\right].$$
Les inclusions $E_{s+\s} \subset E_s$ induisent des applications de restriction.  
\begin{definition} L'approximation définie ci-dessus sera appelée approximation tautologique associée à l'espace vectoriel échelonné $E$.
\end{definition}  
%%%%%%%%%%%%%%%%%%%%%%%%
%%%%%%%%%%%%%%%%%%%%%  

%%%%%%%%%
\subsection{Approximations ultra-violettes}
\label{SS::UV}
On note $C^k_K$ le faisceau des fonctions $k$-fois différentiable au sens de Whitney sur un fermé $K \subset \CM^d$~\cite{Whitney_extension}.
Con\-si\-dé\-rons une famille décroissante de fermés  $(K_{n}),\ n \in \NM $, $K_n \subset \CM^d$ dont on note $K_\infty$ l'intersection. La restriction du faisceau des germes de fonctions holomorphe $\Ot_{\CM^d}$ à chacun des $K_n$ donne lieu à des morphisme de restriction
$$\cdots \to \Ot_{K_n} \to \Ot_{K_{n+1}}  \to \Ot_{K_{n+2}} \to \cdots $$
Pour tout $n$ et pour tout $k$ chacun de ces faisceaux se projette sur celui des fonctions $C^k$ sur $K_\infty$~:
$$ \Ot_{K_n} \to C^k_{K_\infty}$$
pour tout entier $k$.

Fixons un système de voisinage $B_s,\ s \in [0,S]$ croissant, compacts d'un point $x \in K$. Ces voisinages donnent lieu à des structures d'espaces vectoriel échelonné comme en \ref{SS::echelon} sur les anneaux locaux $\Ot_{K_n,0} $, $ C^k_{K_\infty,0} $ et par conséquent à des approximation de ce dernier~:
$$E_{n,s}:=\G(\intr(K_n \cap B_s),\Ot_{\CM^d}) \cap C^k(K_n \cap B_s,\CM) . $$

 %%%%%%%%%%%%%%%%%%%
 %%%%%%%%%%%%%%%%%%%%%%%%%

  %%%%%%%%%%%%%%%%%%%%%%%
\subsection{Produits infinis}
Soit  $E:=(E_n)$ une  approximation d'un espace vectoriel échelonné $E_\infty$. 
\begin{proposition}
\label{P::produits_2}
Soit $u=(u_n)$ un $\tau $-morphisme $1$-borné de $E$ exponentiables.
Si la série numérique
$$\sum_{n \geq 0} N^1_{\tau}(u_n)$$ est convergente  dans $\RM$ alors
la suite $(r\,g_n)$ définie par 
$$g_n=e^{u_n} r_n \cdots r_2 e^{u_1} r_1 e^{u_0}$$
 converge dans  $\Lt(E_0,E_\infty)$.
\end{proposition}
La démonstration est analogue à celle de la proposition~\ref{P::produits}. En voici les grandes lignes.
On fixe $s \in ]0,\tau]$ et on choisit $\l,\mu \in]0,1[$ tels que
 $$ \sup_{n \geq 0}\frac{3}{(1-\l)s}N_{s}^1(u_n) <1$$
 et
 $$ \sup_{n \geq 0}\frac{3}{(1-\mu)\l s}N_{\l s}^1(u_n) <1.$$ 

On  note   $\| \cdot \|_{\l,n}$ la norme d'opérateur dans l'espace $\Lt(E_{s,0},E_{\l s,n})$.
  
 En développant l'exponentielle, on obtient l'estimation~:
 $$\| g_n-r_ng_{n-1}\|_{\l \mu,n} \leq  K_{\l,\mu} \frac{N^1_{\l s}(u_n)}{\l s}$$
 avec 
$$K_{\l,\mu}:=\sup_{n \geq 0} \frac{3 C_\l}{1-\mu-\frac{3N^1_{\l s}(u_n)}{\l s}} $$
et
 $$C_\l:= \prod_{n \geq 0} \frac{1}{1-\frac{3}{(1-\l)s}N_{s}^1(u_n)}.$$

 Cette inégalité montre que pour tout $j \geq 0$, la suite $(rg_n),\ n \geq j$ définit une suite de Cauchy dans l'espace de Banach
 $\Lt(E_{s,0},E_{\l \mu s,\infty})$.  
 %%%%%%%%%%%%%%%%%%%%%%%%%%%%%%%%%%%%%%%%%%%%%
  \section{Le théorème KAM gé\-né\-rali\-sé} 
 %%%%%%%%%%%%%%%%
%%%%%%%%%%%%%%%%%%%%%

%%%%%%%%%%%%%%%%%%%%%%%
 \subsection{\'Enoncé du théorème}
 
 \begin{definition} Une application $f:E \to F$ entre deux  suites d'espaces vectoriel $S$-échelonnés $E=(E_n), F=(F_n)$
 est dite $l$-modérée si la suite
 $$p_n:=\sup\{ \s^l \frac{|f_n(x)|_s}{1+ | x |_{s+\s}} :s \in ]0,S[ ,\s \in ]0,S-s[,x \in E_{n,s+\s} \}$$
 est à croissance modérée. 
 \end{definition}
 Tout morphisme modérée définit bien entendu une application modérée.
\begin{theoreme}
 \label{T::KAM} Soit $E$ une approximation, $a \in E$, $F$ une sous-approximation directe de $E$ et $\alg$ un sous-espace vectoriel de $\Mt^1(E)^{(2)}$. 
 Soit $M$ un sous-ensemble $e^{\alg}$-invariant. Supposons qu'il existe pour certains $k,l \geq 0$ une application  $l$-modérée 
 $$j:F \mapsto \Mt^k(\Gr(M)/\Gr(F),\alg)$$
 telle que   $j(\a)$ soit un quasi-inverse du morphisme
$$ \alg \to \Gr(M)/\Gr(F),  u \mapsto \overline{u (a+\a)} .$$
 Pour tout $x \in a+M _0$, il existe une suite $(u_n)$, avec $u_n \in (\alg_n)^{n-k}$, telle que
 \begin{enumerate}[{\rm i)}]
 \item la suite $g_n:=re^{u_n} r_n \dots e^{u_2} r_2 e^{u_1} r_1  e^{u_0}  $ converge vers un élément $g\in \Lt(E_0,E_\infty)$ ;
 \item $  g(x)=r(a) (\mod F_\infty). $
 \end{enumerate}
 \end{theoreme}
  %%%%%%%%%%%
 \subsection{Principe de la démonstration du théorème \ref{T::KAM}}
 Pour chaque $b \in M$, nous allons construire de proche en proche des suites  $(u_n)$, $(b_n)$, $(\a_n)$, $(c_n)$.
Posons $$a_0=a, b_0=b, u_0=j_0(0)\bar b_0.$$

Fixons un supplémentaire $G$ de $F$ dans $E$.
Soit $\a_0 \in F_0$ et $c_0 \in G_0$ tels que
$$\a_0+c_0=b_0-u_0(a_0).$$

On construit les termes suivants par les formules~:
  \begin{enumerate}[{\rm 1)}]
 \item $a_{n+1}=r_{n+1}(a_n+\a_n)$ ;
 \item $b_{n+1}=r_{n+1}e^{-u_n} (a_n+b_n)-a_{n+1}$ ;
 \item $u_{n+1}=j_{n+1}(\sum_{i=0}^n\a_i)\bar b_{n+1}$ ;
  \item $\a_{n+1}+c_{n+1}=b_{n+1}-u_{n+1} ( a_{n+1}),\ \a_{n+1} \in F_{n+1},\ c_{n+1} \in G_{n+1}$.
 \end{enumerate}

Nous allons voir que  $(g_n)$ converge vers une limite $g$,  $(u_n)$ et $(c_n)$ tendent vers $0$ et $(\sum r\a_n)$ converge vers une limite $\a \in F_\infty$. 

Dans ce cas, la suite $(rg_n x)$ converge vers $x' \in a+F$. En effet, 
$$\lim_{n \to +\infty}rc_n+\lim_{n \to +\infty}r\a_n=0 $$
et par définition de $j$, on a~:
$$ b_n+\a_n+c_n=u_n(a_n).$$
En prenant l'image par $r$ dans les deux membres de l'égalité et en passant à la limite sur $n$,  on trouve 
$$\lim_{n \to +\infty}r(b_n)=0\ (\mod F_\infty).$$
Par ailleurs
$$a_{n+1}+b_{n+1}=r_{n+1}g_n (a+b) $$ 
donc, en prenant, à nouveau, l'image par $r$ et la limite sur $n$, on trouve bien
$$r(a)=g (a+b)\ (\mod F_\infty).$$ 
CQFD.

%%%%%%%%%%%%%%%%%%%%%%%%%%%%%%%%%
 \subsection{Démonstration du Théorème \ref{T::KAM}}

\begin{lemme}
 \label{L::reste}
Soit $E$ un espace vectoriel échelonné.
 Pour tout $\tau$-morphisme $1$-borné $u$, tout $s \in ]0,\tau[$ vérifiant la condition 
 $$\frac{3N^1_\tau(u)}{(\tau-s)} \leq \frac{1}{2} $$
 et tout $x \in E_\tau$, on a les inégalités
 \begin{enumerate}[{\rm 1)}]
 \item $\displaystyle{| (e^{-u}(\Id +u)-\Id)  x|_s \leq   \frac{36|x|_\tau}{(\tau-s)^2} N^1_\tau(u)^2}$ ;
 \item $\displaystyle{ | (e^{-u}(\Id +u)-\Id)  x|_s \leq   \frac{2| u(x)|_\tau}{(\tau-s)} N^1_\tau(u)}$ ;
 \item $\displaystyle{| (e^{-u}-\Id)  x|_s \leq  \frac{6 |x|_\tau}{(\tau-s)} N^1_\tau(u)}$ ;
  \item $\displaystyle{| (e^{-u}-\Id)  x|_s \leq   2  | u(x)|_\tau}$ ;
 \item $\displaystyle{| e^u x |_{s} \leq 2   |x |_\tau}$.
\end{enumerate} 
\end{lemme} 
\`A titre d'exemple, montrons la première de ces inégalités.
L'égalité
$$ e^{-u}(\Id +u)-\Id=\sum_{n \geq 0} \frac{(n+1)}{(n+2)!}(-1)^{n+1}u^{n+2}$$
donne l'estimation
$$ | \sum_{n \geq 0} (-1)^{n+1} \frac{(n+1)}{(n+2)!}u^{n+2} ( x) |_s \leq  | x |_\tau \sum_{n \geq 0}  
\frac{(n+1)3^{n+2}}{(\tau-s)^{n+2}}N^1_\tau(u)^{n+2} .$$
Comme
$$\frac{3N^1_\tau(u)}{\tau-s} \leq 1 ,$$
la série du membre de droite est égale à
$$t^2\sum_{n \geq 0}(n+1)t^n=\frac{t^2}{(1-t)^2},\ {\rm \ avec\ }t= \frac{3N^1_\tau(u)}{\tau-s}.$$
En utilisant l'inégalité
$$\frac{1}{(1-t)^2} \leq 4,\ \forall t \in [0,\frac{1}{2}], $$
on trouve bien la majoration du lemme.
 
Quitte à multiplier toutes les normes par une même constante, on peut supposer que~:
  $$| a |_s \leq \frac{1}{72}. $$  
  Ce qui nous évitera d'avoir à s'occuper des constantes qui inter\-vien\-nent dans les estimations du lemme.
  
    Soit $l$ tel que l'application $j$ soit $l$-modéré.   On définit la suite $(p_n)$ par
$$p_n:=  \sup  \left\{ 2\s^l \frac{N^k_{s}(j_n(\a))}{1+ | \a |_{s+\s}}: \a \in (E_n)_s \right\}$$
si cette borne supérieure est au moins égale à un, si ce n'est pas le cas on prend $p_n:=1$.

Les produit infinis
 $$\b_n:=\prod_{ i > 0  } p_{n+i}^{-2^{-i-1}}=p_{n+1}^{-\frac{1}{2}}\, p_{n+2}^{-\frac{1}{4}}\, p_{n+2}^{-\frac{1}{8}} \cdots$$
 définissent une suite majorée par $1$ car $p_n \geq 1$ et, de plus,
 $$\b_{n+1}=p_{n+1} \b_n^2 .$$
Par ailleurs, comme $(p_n)$ est  à croissance modérée, on a~: $\b_n>0$ pour tout $n \geq 0$.

Définissons à présent les suites $(s_n)$, $(\s_n)$ par $\displaystyle{\s_n=\frac{s}{3^{n+2}}}$ et
$$ s_{n+1}=s_n-3\s_n,  s_0=s, s_1=\frac{2s}{3},\dots.$$

Pour $n>0$, le vecteur  $b_n$ s'écrit sous la forme
$$b_{n+1}=r_{n+1}(A_n+B_n+C_n), $$
avec
$$A_n:=(e^{-u_n}(\Id +u_n)-\Id) a_n,\ B_n:= (e^{-u_n}-\Id)\a_n,\ C_n:=e^{-u_n} c_n.$$

Comme $j$ définit un quasi-inverse, on a~:
$$\ord(u_na_n)=\ord(b_n),\ \ord(\a_n) \geq \ord(b_n),\ \ord(c_n) > \ord(b_n). $$
La deuxième inégalité du lemme montre que
$$\ord(A_n) \geq \ord(u_n(a_n))+\ord(u_n)-1>\ord(b_n).$$
De même
$$ \ord(B_n) \geq \ord(u_n \a_n) \geq \ord(u_n)+\ord(\a_n)>\ord(b_n) $$
 
 Fixons $m \geq 0$ pour que la projection sur $F$ soit $m$-bornée de norme égale à~$1$.
  
Les ordres de $A_n$, $B_n$, $C_n$, $c_n$ et $\a_n$ croissent avec $n$, donc, quitte à remplacer $s$ par $s'<s$ suffisament petit  et  les éléments $a_0$ et $b_0$ par $a_N$ et $b_N$, avec $N$ assez grand, on peut supposer  que les inégalités suivantes sont vérifiées pour $n=0$~:
 \begin{enumerate}[{\rm i)}]
 \item $\displaystyle{| A_n|_{s_n} \leq \frac{\b_n^2 \s_n^{2k+2l+2m+3} }{3}} $ ;
  \item $\displaystyle{| B_n|_{s_n} \leq \frac{\b_n^2 \s_n^{2k+2l+2m+3} }{3}} $ ;
    \item $\displaystyle{| C_n|_{s_n} \leq \frac{\b_n^2 \s_n^{2k+2l+2m+3} }{3}} $ ;
 \item $\displaystyle{| \a_n |_{s_n} \leq \b_{n} \s_{n-1}^{k+l+m+1}} $ ;
    \item $\displaystyle{| c_n |_{s_n} \leq  \b_{n} \s_{n-1}^{k+l+m+1 }}$ ;
\item $\displaystyle{s \leq \frac{1}{3^{2k+2l+4m+2}}  }$.
 \end{enumerate}
 
Montrons par récurrence sur $n$ qu'elles sont alors vérifiées  pour tout $n \geq 0$.

Supposons ces inégalités vérifiées jusqu'à l'ordre $n$, $n \geq 0$. 
Comme $b_{n+1}=r_{n+1} (A_n+B_n+C_n)$, d'après i), ii) et iii) au rang $n$, on a~:
$$|b_{n+1}|_{s_{n+1}} \leq   \b_{n}^2\s_{n}^{2k+2l+3m+3} .$$

 En utilisant iv) jusqu'au rang $n$, on a~:
 $$\sum_{i=0}^n| \a_i |_{s_i} \leq 1 .$$
 Par définition de la suite $(p_n)$, on a~:
$$N_{s_{n+1}+2\s_n}^k (j_{n+1}(\sum_{i=0}^{n}\a_i)) \leq \frac{p_{n+1}}{2\s_{n}^l}(1+\sum_{i=0}^{n}| \a_i |_{s_{n+1}+\s_{n}})  \leq \frac{p_{n+1}}{\s_n^l}  .$$
 
 Comme 
 $$u_{n+1}=j_{n+1}(\bar b_{n+1}), $$ 
 on obtient ainsi l'estimation
 $$(*)\ N^1_{s_{n+1}+\s_n}(u_{n+1}) \leq \frac{p_{n+1}\b_{n}^2}{\s_n^{k+l}}  \s_n^{2k+2l+2m+3}= \b_{n+1}  \s_n^{k+l+2m+3}.$$

Montrons à présent les inégalités iv) et v) au rang $n+1$. L'inégalité 
$$\sum_{i=0}^{n}| \a_i |_{s_i} \leq 1 $$
donne l'estimation
$$ |u_{n+1} ( a_{n+1}) |_{s_{n+1}} \leq \frac{2 N^1_{s_{n+1}+\s_n}(u_{n+1})}{\s_n} .$$
D'après (*), on a alors
$$|u_{n+1} ( a_{n+1}) |_{s_{n+1}}  \leq 2\b_{n+1}\s_n^{k+l+2m+2} . $$
 Comme $\b_n^2 \leq \b_{n+1}$ et $s<1$,  cette estimation et celle sur $b_{n+1}$, nous donne~:
 $$|b_{n+1}-u_{n+1} ( a_{n+1})|_{s_{n+1}} \leq  |b_{n+1}|_{s_{n+1}} +|u_{n+1} ( a_{n+1})|_{s_{n+1}} \leq 3\b_{n+1} \s_n^{k+l+2m+2}.  $$
L'hypothèse vi) entraîne que
$$\s_{n+1} \leq s \leq  \frac{1}{9}, $$
d'où l'estimation
$$|b_{n+1}-u_{n+1} ( a_{n+1})|_{s_{n+1}}   \leq  \frac{\b_{n+1} \s_n^{k+l+2m+1}}{3} .$$
La projection sur $F$ est $m$-bornée de norme 1, on a donc~:
$$|\a_{n+1}|_{s_{n+1}} <   \frac{\b_{n+1} \s_n^{k+l+m+1}}{3}$$
et
$$|c_{n+1}|_{s_{n+1}} \leq |b_{n+1}-u_{n+1} ( a_{n+1})|_{s_{n+1}}+|\a_{n+1}|_{s_{n+1}}\leq \b_{n+1} \s_n^{k+l+m+1}.$$
Ce qui démontre iv) et v) à l'ordre $n+1$. 

 En appliquant le point 1) du lemme avec
$$A_{n+1}:=e^{-u_{n+1}}(\Id +u_{n+1})-\Id) a_{n+1} $$
et $\tau-s=\s_n$, on obtient l'inégalité
$$| A_{n+1}|_{s_{n+1}} \leq  \frac{N^1_{s_{n+1}+\s_n}(u_{n+1})^2}{\s_n^2} .$$
En utilisant l'inégalité (*), on trouve alors
$$| A_{n+1}|_{s_{n+1}} \leq   \frac{\b_{n+1}^2  \s_n^{2k+2l+4m+6}}{\s_n^2} =\b_{n+1}^2  \s_{n}^{2k+2l+4m+4} . $$
Comme
$$ \s_n  \leq \frac{s}{9} \leq \frac{1}{3^{2k+2l+4m+4}} ,$$
on a~:
$$ \s_{n}^{2k+2l+4m+4}  \leq  \frac{1}{3}\s_{n+1}^{2k+2l+4m+3}  . $$
Nous avons donc démontré i) au rang $n+1$.

 Appliquons maintenant le point 2) du lemme à $B_{n+1}$  avec $\tau-s=\s_n$. On obtient ainsi l'inégalité~:
$$| B_{n+1}|_{s_{n+1}} \leq  \frac{ N_{s_{n+1}+\s_n}(u_{n+1})}{\s_n} |\a_{n+1}|_{s_{n+1}+\s_n}.$$
L'estimation iv) au rang n+1 et l'estimation (*) entraînent l'inégalité~:
$$| B_{n+1}|_{s_{n+1}} \leq  \b_{n+1}^2 \s_{n}^{2k+2l+3m+4}   . $$
On a~:
$$ \s_n \leq  \frac{s}{9}  \leq   \frac{1}{3^{2k+2l+3m+4}}  $$
d'où l'inégalité
$$\s_{n}^{2k+2l+3m+4}  \leq\frac{1}{3} \s_{n+1}^{2k+2l+3m+3}. $$
Nous avons donc démontré ii) au rang $n+1$.

Appliquons maintenant le point 3) du lemme à $C_{n+1}$  avec $\tau-s=\s_n$. 
On obtient ainsi l'inégalité~:
$$| C_{n+1}|_{s_{n+1}} \leq  \frac{N_{s_{n+1}+\s_n}(u_{n+1}) }{\s_n} |c_{n+1}|_{s_{n+1}+\s_n}.$$
L'hypothèse de récurrence v) et l'estimation (*) entraînent l'inégalité~:
$$| C_{n+1}|_{s_{n+1}} \leq  \frac{1}{3}\b_{n+1}^2 \s_{n+1}^{2k+2l+3m+3}  . $$
Nous avons ainsi démontré iii) au rang $n+1$, ce qui achève la démonstration du théorème.
     %%%%%%%%%%%%%%%%%%%%%   
   %%%%%%%%%%%%%%
%%%%%%%%%%%%%%%%%%%%%%%%%%%%
\section{Forme normale d'un hamiltonien en un point critique}
\subsection{Compléments sur la forme normale de Birkhoff}
  %%%%%%%%%%%%%%% 
Con\-si\-dé\-rons à nouveau l'espace $\RM^{2n}$ muni de coordonnées $q_i,p_i,\ i=1,\dots,n $ et de  la forme symplectique standard~:
$$\omega:=\sum_{i=1}^n dq_i \w dp_i.$$
  Supposons les $\a_i$ linéairement indépendants sur $\QM$. Dans ce cas, pour tout $l$,  il existe un symplectomorphisme 
 $$\p_l:(\RM^{2n},0) \to (\RM^{2n},0) $$
  et un polynôme 
 $$A_l(X_1,X_2,\dots,X_n) \in \RM[X_1,\dots,X_n], $$ de degré $l$, appelé {\em polynôme de Birkhoff}, tels que~:
$$H \circ \p_l = A_l(p_1^2+q_1^2,\dots,p_n^2+q_n^2) +o(2l) ,\ \a_i \in \RM .$$
 
En prenant la limite sur $l$, on obtient des séries formelles $A,\p$ qui vérifient
$$ H \circ \p = A(p_1^2+q_1^2,\dots,p_n^2+q_n^2) .$$
La série $A$ est unique. On l'appelle la {\em forme normale de Birkhoff}.  

 Con\-si\-dé\-rons {\em les applications des fréquences} 
$$\del A_l=(\d_{X_1}A_l,\dots,\d_{X_n}A_l):\RM^n\to \RM^n.$$
 
\begin{definition} L'espace des fréquences de $H$, noté $V(H) \subset \RM^n$, est le plus petit sous-espace affine de $\RM^n$ qui contient les images de $\del A_l $ pour tout $l \geq 2$. 
\end{definition}
Dans les cas gé\-né\-riques,  les applications des fréquences sont des isomorphismes et on a~: $V(H)=\RM^n $. Ce sont les cas que l'on appelle {\em isochroniquement non dégé\-né\-rés}. 

Soit $e_1,\dots ,e_d \in \RM^n$ une base de l'espace des fréquences. 
 Nous allons re-écrire les fonctions $A_l$ en faisant intervenir directement les applications des fréquences. Pour cela,
munissons l'espace produit
$$\RM^{2n} \times \RM^d=\{(q,p,\l) \}$$  
de la structure de Poisson induite par le bivecteur
$$v=\sum_{i=1}^n \d_{q_i} \w \d_{p_i} .$$ 
Considérons l'application
$$f:(\RM^{2n},0)  \to (\RM^n,0),\ (q,p) \mapsto (p_1^2+q_1^2, p_2^2+q_2^2,\dots, p_n^2+q_n^2)$$
et notons $I$ l'idéal de $\RM[[q,p,\l]]$ engendré par les composantes de l'application
$f-\sum_{i=1}^d\l_i e_i .$

Si deux fonctions sont égales modulo le carré de l'idéal $I$ alors elles définissent le même flot hamiltonien sur les variétés
$$L_\l=\{ (q,p) \in \RM^{2n}:f=\sum_{i=1}^d \l_i e_i \}. $$
En effet, pour $g,h \in I$, on a~:
$$\{ H+gh,-\}=\{ H,-\}+g\{h,-\}+h\{g ,-\}=\{ H,-\}\ (\mod I) .$$

La formule de Taylor à l'ordre $2$ donne l'égalité~(lemme d'Hadamard)~:
$$A(f)=A(x)+(\del A(x),f-x)  \ (\mod I^2),\ x=\sum_{i=1}^d \l_i e_i. $$
Par conséquent, les fonctions $A(f) $ et $(\del A,f)$ définissent la même dérivation hamiltonienne de l'anneau $\RM[[q,p,\l]]$.

On peut ainsi ré-écrire la normalisation de Birkhoff sous la forme
$$\p(H)=\sum_{i=1}^n (\d_{X_i} A(x))p_iq_i \ (\mod I^2+\RM[[\l]]) $$
où $\p$ est un automorphisme de $\RM[[\l,q,p]]$ qui préserve la structure de Poisson. 

La signature de la partie quadratique de la forme normale de Birkhoff ne joue pas de rôle particulier. On peut considérer une fonction holomorphe
$$H:(\CM^{2n},0) \to (\CM,0)$$
avec un point critique de Morse à l'origine. Il existe alors des séries formelles $A,\p$ qui vérifient
$$ H \circ \p = A(p_1q_1,\dots,p_nq_n) .$$
Si $H$ est réelle, la partie quadratique de sa forme norme de Birkhoff est conjuguée par une application linéaire à la partie quadratique de $A$.
Cette forme linéaire envoient la conjugaison complexe sur une involution antiholomorphe $\tau$. Les applications $A,\p$ envoie $\tau$
 sur la conjugaison complexe. La notion d'espace des fréquences se définit dans le cas complexe comme dans le cas réel  et lorsque $H$ est réel, cet espace vectoriel est la complexification de celui définit  sur le corps des nombres réels. 
 
 Par la suite, nous ne préciserons pas ces structures réelles, car elles ne jouent aucun rôle dans la démonstration.

  %%%%%%%%%%%%%%%%%%%%%%%%% 
  \subsection{Produits de Hadamard}
\label{SS::Hadamard} 
 Con\-si\-dé\-rons l'espace $\Ot_{\CM^d,0}$ munit de l'échelonnement en espace de Hilbert décrit au n°\ref{SS::echelon}. Con\-si\-dé\-rons
 l'approximation tautologique associée à cet échelonnement et posons 
 $$z^i=z_1^{i_1}z_2^{i_2}\dots z_n^{i_n},\ i \in \NM^d.$$ 
L'idéal maximal  $\Mt_{\CM^d,0} \subset \Ot_{\CM^d,0}$ est ainsi muni d'un échelonnement.
  Le produit de Hadamard de deux séries 
  $$f:=\sum_i a_i z^i,\ g:=\sum_i b_i z^i$$  
  est défini par
  $$f \star g:=\sum_{i \in \NM^n} a_ib_i z^i $$

 \begin{proposition}
 Soit $\a \in \RM^n$ un vecteur qui satisfait la condition de Bruno. Posons 
 $$f:= \sum_{i \in \NM^n} (\a,i) z^i,\ g_n:=\sum_{\ \| i \| \leq n} (\a,i)^{-1}z^i$$ alors 
 les produits de Hadamard
$$g_n \star:  (\Mt_{\CM^d,0})_n \to   F_n,\  \sum_{i \in \NM^n} f_i z^i \mapsto \sum_{\ \| i \| \leq n} \frac{f_i}{ (\a,i)} z^i,\ f_i \in \CM,$$
avec $F=\Mt_{\CM^d,0}$ ou $F=\Gr(\Mt_{\CM^d,0})$,
 définissent un quasi-inverse $0$-modéré  du produit de Hadamard par $f$ dont la norme est majoré par la suite $(\s(\a)_n^{-1})$.
\end{proposition}
L'application définit évidemment un quasi-inverse.
Par ailleurs, on a~:
$$| u_n(f)|_{s}^2 = \sum_{\ \| i \| \leq n} \frac{|f_iz^i|^2_{s}}{ |(\a,i)|^2} \leq \frac{1}{[\a]_n^2}\sum_{\ \| i \| \leq n} |f_i|^2|z^i|^2_{s}$$
ce qui donne l'estimation souhaitée pour la norme de cet inverse.

Soit à présent $K$ la famille de compact
$$K_m =  (\Dt_a)_m,\ K_\infty =\Dt_a $$
et $B_s,\ s \in ]0,S]$ un système fondamental de voisinages croissants compacts de $\a$.
Munissons l'espace $C^k_{K_\infty,\a} $ de l'approximation
$$(C^k_{K_\infty,\a})_{n,s}:=C^k(K_n \cap B_s,\CM). $$
\begin{proposition}
 \label{P::Hadamard}
 Si $a=(a_k)$ est une suite à décroissance modérée alors les applications définies par
$$u_n: (C^k_{K_\infty,\a}  \hat \otimes \Mt_{\CM^d,0})_n  \to  F_n  ,\ \sum_{i \in \NM^n} f_i \otimes z^i \mapsto \sum_{\| i \| \leq n} \frac{f_i}{ (\a,i)} \otimes z^i ,$$
avec $F=C^k_{K_\infty,\a}  \hat \otimes \Mt_{\CM^d,0}$ ou $F=\Gr(C^k_{K_\infty,\a}  \hat \otimes \Mt_{\CM^d,0})$,
 définissent un quasi-inverse  $k$-modéré de
 $$v: \sum_{i \in \NM^n} f_i \otimes z^i \mapsto \sum_{i \in \NM^n} (\a,i)  f_i \otimes z^i,$$
 dont la norme est au plus $(a_k^{-1})$.
\end{proposition}
Si $k=0$ alors la démonstration est identique à celle du cas précédent. En effet~:
$$| u_n(f)|_{s}^2 = \sum_{\ \| i \| \leq n} \frac{|f_i \otimes z^i|^2_{s}}{ |(\a,i)|^2} \leq \frac{1}{a_n^2}\sum_{\ \| i \| \leq n} |f_i|^2| z^i|^2_{s}$$
ce qui donne l'estimation souhaitée pour la norme de $u_n$. 
 
Fixons $n$ et notons provisoirement $\|\cdot \|_{k,s}$ la norme de $ (C^k_{K_\infty,\a}  \hat \otimes \Mt_{\CM^d,0})_{n,s}$.
En vertu des inégalités de Cauchy, on a~:
$$ \| f \|_{k,s} \leq \frac{k!}{\s^k}\| f \|_{0,s+\s} ,\ \s <1,$$
ce qui implique la propriété pour $k >0$.

L'injection $(H_s) \to (E_s)$ d'espace vectoriel échelonnée est $1$-bornée par conséquent les deux propositions précédentes restent valables pour la structure échelonnée $(E_s)$ à condition d'augmenter l'indice de mo\-dé\-ration de $0$ à $1$ pour la première et de $k$ à $k+1$ pour la seconde.

%%%%%%%%%%%%%%%%%%%
\subsection{Forme normale sur la fibre spéciale}
Notons $I \subset \Ot_{\CM^{2n},0}$ l'idéal engendré par les germes de fonction
$$p_1q_1,\ p_2q_2,\dots,\ p_nq_n.$$
   \begin{proposition} 
   \label{P::fibre}
Soit $H:(\CM^{2n},0) \to (\CM^n,0)$ un germe de fonction analytique de la forme
 $$H(q,p)=\sum_{i=1}^n \a_i p_iq_i+o(2).$$
   Si le vecteur $\a=(\a_1,\dots,\a_n)$ satisfait la condition de Bruno alors pour toute fonction de la forme
 $$H+R,\ R \in \Mt^3_{\CM^{2n},0} $$
 il existe un automorphisme symplectique $\p \in \Aut(\Ot_{\CM^{2n},0})$ tel que
 $$\p(H+R)=H \ (\mod I^2). $$
 De plus, si $H$ et $R$ sont invariantes par une involution réelle alors $\p$ peut être choisit réel.
 \end{proposition}
En particulier, dans un voisinage suffisamment petit de l'origine, $H+R$ admet une variété lagrangienne complexe invariante symplectomorphe à
$$\{ (q,p): p_1q_1=p_2q_2=\dots=p_nq_n=0 \}. $$
 
 Cette proposition est une application directe du théorème KAM gé\-né\-ralisé.
 En voici les détails, munissons l'espace vectoriel
 $$E:=\Ot_{\CM^{2n},0}$$
 de la structure échelonnée $(H_s)$ définie en \ref{SS::echelon}. Celle-ci induit une approximation tautologique sur $E$. 
Désignons  par $\Mt_{\CM^{2n},0}$ l'idéal maximal de l'anneau local $\Ot_{\CM^{2n},0} $.
 Munissons les espaces vectoriels
 $$ M:=\Mt^3_{\CM^{2n},0},\ F:=I^2 \cap \Mt^3_{\CM^{2n},0} $$
 de la structure échelonnée induite par $E$ ainsi que de leurs approximations  tautologiques.
 L'espace vectoriel $\Gr(M)$ (resp. $ \Gr(F)$) s'identifie aux sous-espace des polynômes dans les variables $q,p$ contenus dans $M$ (resp. dans $F$).
 
 Le groupe $G$ des germes de symplectomorphismes dont la partie linéaire est l'identité agit sur $H+\Mt^3_{\CM^{2n},0}$.
 Désignons par $\alg$ les dérivations de la forme
$$\{  h,\cdot \} \in \Der(\Ot_{\CM^{2n},0}),\ h \in   \Mt^3_{\CM^{2n},0}. $$
 
On construit un inverse $j(\a)$ modéré de
$$\rho(\a):\alg \to M/F \mapsto   g \mapsto \{ g,H+\a \},\ \a \in I^2$$
en cherchant d'abord un inverse modulo $I$ puis en ajoutant une correction afin  d'obtenir l'inverse modulo $I^2$. 

Pour cela, commençons par remarquer que la base $p_1q_1,\dots,p_nq_n$ permet de scinder la suite exacte d'espaces vectoriels
$$0 \to I/I^2 \to \Ot_{\CM^{2n},0}/I^2 \to \Ot_{\CM^{2n},0}/I \to 0 .$$
On obtient ainsi un isomorphisme d'algèbres
$$ \Ot_{\CM^{2n},0}/I^2 \approx \Ot_{\CM^{2n},0}/I \oplus I/I^2   $$
et par suite un isomorphisme d'espaces vectoriels
$$M/F \approx M/(I \cap M) \oplus (I \cap M)/(I^2 \cap M). $$
Dans la suite, nous identifierons chaque espace quotient de la forme $A/B$ à l'orthogonal de $B$ dans $A$.

 Soit, à présent, $g_n$ les fonction définies par~:
$$g_n =\sum_{\| i \|  \leq n,\ i \neq 0} \frac{1}{(\l,i)} (pq)^i,\ (pq)^i:=(p_1q_1)^{i_1}(p_2q_2)^{i_2} \dots (p_nq_n)^{i_n}.$$
Les résultats du \ref{SS::Hadamard}  montrent que les produit de Hadamard par les  $g_n$ définissent un quasi-inverse à droite $0$-modéré de
$$\Mt_{\CM^{2n},0}/I \to \Mt_{\CM^{2n},0}/I,\ a \mapsto \{ a,H \}.  $$
Cette propriété permet de définir un quasi-inverse à droite $j(\a)$ de $\rho(\a)$ par la formule~:
$$  M_n/(I \cap M_n) \oplus (I \cap M_n)/(I^2 \cap M_n) \to \alg_n:(a,b) \mapsto a \star g_n+b \star g_n-\{ a \star g_n ,\a \} \star g_n. $$
En effet~:
$$\{ a \star g_n+b \star g_n-\{ a \star g_n ,\a \} \star g_n,H+\a \}=a+b+\{ a \star g_n,\a \}-\{\{ a \star g_n ,\a \} \star g_n,H \} \ (\mod I^2). $$
et la somme des deux derniers termes est nulle, car~:
$$\{\{ a \star g_n ,\a \} \star g_n,H \} =\{ a \star g_n ,\a \}  .$$
Ceci montre que $j(\a)$ est un quasi-inverse à droite $0$-modéré de $\rho(\a)$ associée à la structure échelonné $(H_s)$ donc $1$-modéré pour la structure $(E_s)$.

L'application $\a \mapsto j(\a)$ ne fait intervenir que des dérivées du premier ordre en $\a$. D'après les inégalités de Cauchy, elle est donc  $1$-modérée. Les conditions du théorème KAM gé\-né\-ralisé sont satisfaites, la proposition est démontrée.

%%%%%%%%%%%%%
\subsection{Forme prénormale}
Notons  
$$\l=(\l_1,\dots,\l_d),\ \mu=(\mu_1,\dots,\mu_d)$$ les coordonnées sur l'espace vectoriel $\CM^{2d}$.
Le bivecteur
$$v=\sum_{i=1}^n \d_{q_i} \w \d_{p_i} $$ 
induit une structure de Poisson sur $\CM^{2n+2d}=\{ (q,p,\l,\mu) \}$.

Soit $e_1,\dots,e_d \subset \CM^n$ une base de $V(H)$.
Notons $I_{\l,\mu} \subset \Ot_{\CM^{2n+2d},0}$ l'idéal engendré par les composantes de l'application
$$g:=(p_1q_1,\dots,p_nq_n)-\sum_{i=1}^d (\l_i-\mu_i)e_i,\ i=1,\dots,n$$
et $I_\mu$ pour l'idéal    de $\Ot_{\CM^{2n+d},0}$ engendré par les restrictions des composantes de $g$ à $\l=0$.

 La projection sur les coordonnées $(q,p)$ (resp. $(q,p,\l)$) induit un morphisme d'anneaux
$$\Ot_{\CM^{2n},0} \subset \Ot_{\CM^{2n+2d},0},\ {\rm\ (resp.\  }\Ot_{\CM^{2n+d},0} \subset \Ot_{\CM^{2n+2d},0} { \rm ) }$$
qui nous permet d'identifier $\Ot_{\CM^{2n},0}$  (resp. $\Ot_{\CM^{2n+d},0}$) avec son image dans $\Ot_{\CM^{2n+2d},0}$.

Par ailleurs, pour tout $k$  il existe un germe de symplectomorphisme 
 $$\p_k:(\CM^{2n},0) \to (\CM^{2n},0) $$
et  un polynôme $A_k \in \RM[X_1,X_2,\dots,X_n]$ de degré $k$  tels que~:
$$H \circ \p_k = A_k(p_1q_1,p_2q_2,\dots,p_nq_n)+o(2k).$$

La proposition \ref{P::fibre} permet donc, sans perte de gé\-né\-ralité, de supposer que
$$\left\{ \begin{matrix}  H  &=& \sum_{i=1}^n \a_i p_iq_i+R,\ R \in I^2\ ;\\ 
  H & = &  A_k(p_1q_1,p_2q_2,\dots,p_nq_n)+o(2k)    \end{matrix} \right.  $$
où $k$ est choisit suffisamment grand pour que l'image de $\del A_k$ soit égale à $V(H)$.
 
%%%%%%%%%%%%%%%%%%%%%%%%%%%%%%%%%%%%
\subsection{Fin de la démonstration}
Posons 
$$f=(p_1q_1,p_2q_2,\dots,p_nq_n) $$
et notons $(\cdot,\cdot)$ le produit scalaire.
Définissons la fonction
$$G= \sum_{i=1}^n \a_i p_iq_i +\sum_{i=1}^d \mu_i (f,e_i)+R(q,p).$$
La forme prénormale de $H$ montre que la restriction  de $G$ à $\l=0$ est égale à~:
$$G_0:= \sum_{i=1}^n \a_i g_i+ \sum_{i=1}^d \mu_i (f,e_i)  \ (\mod I_\mu^2)$$

Nous allons à présent appliquer le théorème KAM gé\-né\-rali\-sé.

Posons $$ L_\infty:=\CM^n \times (\Dt_a-\a).$$
et  munissons l'espace vectoriel $C^{2k}_{L_\infty,0}$ de l'approximation  ultraviolette associée aux ensembles fermés~(voir \ref{SS::UV})~:
$$L_m:=\CM^n \times ((\Dt_a)_m-\a), $$
et aux voisinages
$$B_s=\{ (\l,\mu):\| (\l,\mu)\| \leq s^2 \}. $$ 
Le produit tensoriel $E_\infty:=C^{2k}_{L_\infty,0} \hat \otimes \Ot_{\CM^{2n},0}$ est ainsi muni d'une approximation $(E_m)$ 
Soit 
$$ M_\infty \subset (C^{2k}_{L_\infty,0} \hat \otimes \Ot_{\CM^{2n},0})^{(5)} $$
le sous-ensemble définit par
$$R \in M_\infty \iff V(H+R)=V(H). $$

Posons
$$  F_\infty:=(C^{2k}_{L_\infty,0} \hat \otimes I_{\l,\mu}^2+C^{2k}_{L_\infty,0})^{(5)} $$
et prenons pour $\alg_\infty$ l'espace des dérivations de la forme
$$\{ f,- \}+\sum_{i=1}^k a_i \d_{\mu_i} $$
avec $f \in M_\infty$, $a_i \in (C^{2k}_{L_\infty,0})^{(3)} $.
 
 L'approximation $(E_m)$ induit des approximation sur $M_\infty$ et $F_\infty$.
Décomposons les $\Gr(M_m)/\Gr(F_m)$ dans une somme d'espaces vectoriels  $N_m$, $P_m$. Les sous-espaces vectoriels $N_m$ sont engendrés par les classes des $a(\l,\mu) p^iq^j,\ i \neq j $ et $P_m$ par les classes des $a(\l,\mu)(f,e_i) $. 

Soit
$$g_m =\sum_{\| i \|  \leq m, i \neq 0} \frac{1}{(\l,i)} (qp)^i \in E_m.$$
La décomposition
$$\Ot_{\CM^{2n},0}/I^2 \approx \Ot_{\CM^{2n},0}/I \oplus I/I^2  $$ 
permet d'écrire le quotient $N/(\Gr(F) \cap N)$   comme une somme.
On définit comme précédemment les morphismes
$$A_m(\a): (a,b) \mapsto a \star g_m+b \star g_m-\{ a \star g_m ,\a \} \star g_m. $$

Définissons à présent l'application 
$$B_m:P_m \to \alg_m$$ en prenant pour image de $a(\l,\mu)(f,e_i)$ la dérivation
 $$a(\l,\mu)\d_{\mu_i} \in \alg_m.$$ Les applications 
$$j_m(\a):=A_m(\a)+B_m$$  
définissent  des quasi-inverses à droite $2k$-modérés de
$$\alg \to E,\ g \mapsto \{ g,H+\a \} $$
Par conséquent, d'après le théorème KAM gé\-né\-ralisé, il existe une suite de morphismes de Poisson
$$(\p_m),\ \p_m \in \Lt(E_0,E_m) $$
dont la restriction converge vers un morphisme $\p_\infty$ tel que
 $$\p_\infty(G)=G_0. $$
 
  D'après le théorème d'extension de Whitney, ce morphisme est obtenu par restriction à $(\Dt_a -\a) \times \CM^{2n+d} $ d'une application de classe $C^k$~\cite{Whitney_extension}:
$$\p:(\CM^{2d} \times \CM^{2n},0) \to  (\CM^{2d} \times \CM^{2n},0), (\mu,\l,q,p) \mapsto (\p_1(\l,\mu),\p_2(\l,\mu,q,p)).$$
 Cette application est tangente à l'identité donc d'après le théorème des fonctions implicites, dans un voisinage de l'origine, l'image inverse du sous espace vectoriel
 $$ \{ \mu=0 \} \subset \CM^d $$
est le graphe d'une fonction $C^{2k}$~:
 $$(a_1,\dots,a_d):(\CM^d,0) \to (\CM^d,0). $$
 
La série de Taylor à l'origine de l'application
$$a:(\CM^d,0) \to V(H) \subset \CM^n,\ \l \mapsto \a+\sum_{i=1}^d a_i(\l)e_i $$ 
est de la forme $\del A_k +o(2k)$. Les dérivées en l'origine d'ordre $\leq k$ de $A_k$  engendrent $V(H)$ donc l'application $a$ est non-dégé\-né\-rée. Notons $K$ la préimage de $\Dt_a$ dans un voisinage suffisamment petit de $\a$.
 
 La fonction $a$ est limite uniforme pour la topologie $C^{2k}$ sur $K$ des  $\del A_k$. 
Donc d'après le théorème du graphe fermé, l'application $a$ est le gradient d'une fonction $C^{2k}$~:
$$A:(\CM^n,0) \to (\CM,0) .$$
  
Notons $X$, l'image inverse de $K$ par l'application $f$. On a bien
$$(H \circ \p)_{\mid X}=(G(-,\mu=0) \circ \p)_{\mid X}=(G_0)_{| \mu=\del A(K)}\ (\mod I_{\l,\mu}^2+\CM\{\mu\})=
A(q_1,p_1,\dots,p_nq_n)_{|X}. $$
Ceci achève la démonstration du théorème.
%%%%%%%%%%%%
\subsection{Un théorème de tores invariants}
Comme je l'ai indiqué dans l'introduction, on a des variantes du théorème de forme normale au voisinage d'un tore invariant, ainsi qu'au voisinage d'une orbite périodique~(ce qui implique la conjecture de Herman pour les symplectomorphismes). 

Il reste toutefois encore un énoncé à ajouter à la liste de ces résultats. Notons $M$ le produit de $T^*(S^1)^n \times \RM$
muni de coordonnées «actions-angles»~:
$$\theta_j \in S^1,\ I_j \in \RM, t \in \RM,\ j=1,\dots,n $$
et de la structure de Poisson induite par celle du fibré cotangent au tore~:
$$\sum_{j=1}^n \d_{I_j} \w \d_{\theta_j}. $$ 
Notons $0_M$ le produit de la section nulle du fibré cotangent par le singleton $\{ t=0 \}$ et $\pi$ la fonction
$$\pi:M \to \RM^{n+1} (t,\p,I) \mapsto (t,I). $$
\begin{theoreme} 
Soit $b=(b_i)$ une suite numérique à décroissance modérée et 
$$H:(M,0_M) \to \RM $$
un germe de fonction le long de $0_M$ de la forme
$$H=\sum_{j=1}^n \a_j I_j+ S(I)+ t R(I,\p),\ dS(0)=0$$ avec
$$\a=(\a_1,\dots,\a_n) \in \Dt_b .$$ 
Pour tout $k \geq 0$, il existe une application   $A:(\RM \times \RM^n,0) \to \RM$ de classe $C^k$ et un germe de symplectomorphisme
de classe $C^k$~:
$$\p:(X,0) \to (X',0),\ X:=\pi^{-1} \circ (\del A)^{-1}(\Dt_b),\ \del A:=(\d_{I_1}A,\dots,\d_{I_n} A)$$
tels que
\begin{enumerate}[{\rm i)}]
\item $H \circ \p=A(t,I_1,\dots,I_n) $ ;
\item l'application $\del A:(\RM^{n+1},0) \to (\RM^n,0)$ est non-dégé\-né\-rée ;
\item la restriction de $\p$ aux fibres de $\pi$ est analytique.
\end{enumerate}
\end{theoreme}
 Ce résultat donne une variante théorème KAM classique sans aucune hypothèse de non-dégé\-né\-rescence. On pourra comparer ce résultat avec~\cite{Russmann_KAM} et \cite[Section 2]{Sevryuk_KAM}.

%%%%%%%%
 
%%%%%%%%%%%%%%%%%%%%%%
 \bibliographystyle{amsplain}
\bibliography{master}
 \end{document}